\documentclass{amsart}
\usepackage{amssymb}
\usepackage{verbatim}
\newtheorem{theorem}{Theorem}[section]
\newtheorem{lemma}[theorem]{Lemma}
\newtheorem{proposition}[theorem]{Proposition}

\theoremstyle{definition}
\newtheorem{definition}[theorem]{Definition}

\theoremstyle{remark}
\newtheorem{remark}[theorem]{Remark}

\numberwithin{equation}{section}

\newcommand{\LP}[2]{L^{#1}(X_{#2}, \Sigma_{#2}, m_{#2})}

\newcommand{\M}{\mathcal{M}}

\newcommand{\I}{1\!{\mathrm l}}
\newcommand{\cM}[2]{{\mathcal M}_{#1}\rtimes_{\sigma^{#2}}{\mathbb R}}
\newcommand{\HLP}[2]{L^{#1}(\mathcal{M}_{#2})}
\newcommand{\HGP}[3]{L^{#1}_{#2}(\mathcal{M}_{#3})}
\newcommand{\supp}{{\rm supp}}
\newcommand{\tr}{{\rm tr}}

\begin{document}
\title[Composition Operators on Haagerup $L^p$-spaces]{Composition Operators on
Haagerup $L^p$-spaces}
\author{S Goldstein}
\address{Faculty of Mathematics and Computer Science, {\L}{\'o}d{\'z} University, ul. Banacha 22,
90-238 {\L}{\'o}d{\'z}, Poland} \email{goldstei@math.uni.lodz.pl}
\author{L E Labuschagne}
\address{Department of Mathematical Sciences, University of
         South Africa, P.O.Box 392, 0003 Pretoria, South Africa}
\email{labusle@unisa.ac.za}
\subjclass[2000]{46L52, 47B33}

\date{\today}

\thanks{Goldstein acknowledges the support of the MNiI grant 2 P03A
030 24. Labuschagne was partially supported by a grant under
the joint Poland - South Africa cooperation agreement.}

\keywords{}

\begin{abstract}
Building on the ideas in \cite{L1} we indicate how the concept of
a composition operator may be extended to the context of Haagerup
$L^p$-spaces.
\end{abstract}

\maketitle

\section{Introduction}

Classically a (generalised) composition operator $C$ is a bounded linear
operator $C: \LP{p}{1} \rightarrow \LP{q}{2}$ which in a canonical way is
induced by a non-singular measurable transformation $T : Y \subset X_2
\rightarrow X_1$ from a measurable subset $Y$ of $X_2$ into $X_1$ in the sense
that $C(f)(t) = f \circ T(t)$ if $t \in Y$ and $C(f)(t) = 0$ otherwise. In the
setting of standard Borel spaces, up to sets of measure zero, such non-singular
measurable transformations are in 1-1 correspondence with $*$-homomorphisms
$\LP{\infty}{1} \rightarrow \LP{\infty}{2}$. (See for example the discussion in
section 2.1 of \cite{SM}.) So in the noncommutative world the study of
composition operators on $L^p$-spaces translates to a description and study of
those Jordan $*$-morphisms $J: \M_1 \rightarrow \M_2$ which in some canonical
sense induce a bounded operator $C_J: \HLP{p}{1} \rightarrow \HLP{q}{2}$, where
$\HLP{p}{1}$ and $\HLP{q}{2}$ are the corresponding noncommutative spaces (The
definitions we use will be given in the next section). Now even in the
commutative setting the case $p < q$ tends to be pathological (see
\cite[Corollary, Lemma 1.5]{TY}). In the noncommutative setting one has a
negative result of Junge and Sherman \cite[Corollary 2.7]{JS}. Thus we will
focus on the case where $\infty \geq p \geq q \geq 1$.

At the outset of any self-respecting theory of composition operators
two questions need to be answered: Firstly the question of which
point transformations actually induce composition operators, and
secondly the question of how in the class of all bounded linear maps
from $L^p$ to $L^q$ we may recognise those that come from point
transformations. In our noncommutative endeavour this translates to
firstly identifying those Jordan $*$-morphisms $J: \M_1 \rightarrow
\M_2$ that canonically induce bounded maps $C_J: \HLP{p}{1}
\rightarrow \HLP{q}{2}$, and secondly describing those bounded maps
between noncommutative $L^p$-spaces that come from Jordan
$*$-morphisms. In section 3 we will indicate how the classical
process for constructing composition operators on $L^p$-spaces may
be extended to the setting of von Neumann algebras as well as
indicating a possible answer to the above two questions.

We tried to make the exposition accessible to both specialists in
operator algebras, and also specialists dealing with composition operators
on classical function spaces. This means that in many places we explain
more than is strictly necessary, especially for specialists in operator
algebras. However we do this consciously for the sake of reaching a larger
audience.

We would like to thank David Sherman, who directed our attention
to the paper of Junge and Sherman \cite{JS}, and to the fact that
their Theorem 2.5 on the general form of the (right) $\M$-module
homomorphisms of noncommutative $L^p$ spaces implies our change of
weight result (see Step II in Section 3). It turned out that after
a slight modification we were able to prove their theorem using
our method, at least in the case when $1\leq q\leq p\leq \infty$.
We decided to show the proof to the reader, as it differs
substantially from the proof of Junge and Sherman in that it uses
essentially only duality arguments.

\section{Prerequisites}
Throughout this paper we will assume that $\M_1$ and $\M_2$ are von
Neumann algebras with faithful normal semifinite (\emph{fns} for
short) weights $\varphi_1$ and $\varphi_2$ respectively. For a von
Neumann algebra $\M$ with an \emph{fns} weight $\varphi$, the
crossed product of $\M$ with the modular action induced by $\varphi$
will be denoted by $\cM{}{}$ and the canonical trace on $\cM{}{}$ by
$\tau$. The Haagerup $L^p$ space constructed by means of the action
of $\varphi$ is denoted by $\HGP{p}{\varphi}{}$. Now let $h =
\frac{d\widetilde{\varphi}}{d\tau}$ where $\widetilde{\varphi}$ is
the dual weight on the crossed product. Then $h$ is a closed densely
defined positive non-singular operator affiliated with the crossed
product. In general, $h$ is not $\tau$-measurable, so it has to be
manipulated with caution.

Define, for $q\in[2,\infty[$,
$$ \mathfrak{n}_\varphi^{(q)}=\mathfrak{n}^{(q)} := \{a\in \M: ah^{1/q}
\text{ is closable and } [ah^{1/q}]\in \HGP{p}{\varphi}{}\}$$
(here $[\cdot]$ is used to denote the minimal closure of a given closable operator).
For $p\in[1,\infty[$, denote by $\mathfrak{m}^{(p)}$ the linear span
of elements of the form $b^*a$ with $a,b\in\mathfrak{n}^{(2p)}$.
Then $\mathfrak{n}_\varphi=\mathfrak{n}^{(2)}$, where
$\mathfrak{n}_\varphi=\{a\in\M:\varphi(a^*a)<\infty\}$, and
$\mathfrak{m}_\varphi=\mathfrak{m}^{(1)}\subset \mathfrak{m}^{(p)}$
for each $p>1$, where $\mathfrak{m}_\varphi$ is linearly spanned by
positive elements $a$ from the algebra satisfying
$\varphi(a)<\infty$. The linear extension to $\mathfrak{m}^{(p)}$ of
the map
$$a\mapsto h^{1/(2p)}a^{1/2}\cdot
[a^{1/2}h^{1/(2p)}]:\mathfrak{m}^{(p)}_+\to\HGP{p}{\varphi}{}
$$
is denoted by $\mathfrak{i}^{(p)}$, and the image of $a$ under the mapping by
$h^{1/(2p)}ah^{1/(2p)}$. Other than that, we use the following convention:
whenever a formula consists of (pre)measurable operators only, their
juxtaposition denotes their strong product; otherwise, it denotes the usual
operator product, and we use square brackets for the closure of a closable
operator. Sometimes we add parentheses to avoid ambiguity. For example, if $h$
is not measurable, but $a,b$ and $h^{1/p}b$ are, we write $a(h^{1/p}b)$ to
denote the strong product of $a$ and $h^{1/p}b$.

Let now $X_0$ denote the completion of $\mathfrak{m}_\varphi$ equipped with
the norm $\|a\|_0$ equal to the maximum of $\|a\|$ and
$\|\mathfrak{i}^{(1)}(a)\|_1$. The mappings $\mathfrak{i}^{(p)}$ can be extended
to bounded maps from $X_0$ into $\HGP{p}{\varphi}{}$. Denote by
$\kappa_{p}, 1<p\leq\infty,$ the Banach space adjoint of
$\mathfrak{i}^{(p^*)}$, where $p^*$ is the conjugate index of $p$. Define
additionally $\kappa_{1}$: if $h_\psi$ is an element of
$\HGP{1}{\varphi}{}$ corresponding to the functional $\psi\in\M_*$
(i.e. $h_\psi=\frac{d\widetilde{\psi}}{d\tau})$, then
$\kappa_{1}(h_\psi)$ is an element of $X_0^*$ which maps $a\in
\mathfrak{m}_\varphi$ onto $\psi(a)$. Then $\kappa_{p}$ maps
$\HGP{p}{\varphi}{}$ boundedly into $X_0^*$. If we denote by $X_1$
the closure of
$\kappa_{\infty}(\HGP{\infty}{\varphi}{})+\kappa_{1}(\HGP{1}{\varphi}{})$
in $X_0^*$ and by $L^p(\M,\varphi)$ the image
$\kappa_{p}(\HGP{p}{\varphi}{})$ equipped with the norm $\|a\|_p^\varphi =
\|(\kappa_p^{(\varphi)})^{-1}a\|_p$, then
$L^p(\M,\varphi)=C_{1/p}(X_0,X_1)$, where $C_\theta, 0\leq\theta\leq
1$ is the $\theta$'s interpolation functor for the complex
interpolation method of Calderon. The spaces $L^p(\M,\varphi)$ are
the Terp interpolation spaces. (For a precise explanation of the
interpolation method the reader is directed to Terp's paper
\cite{Tp2}.)

The theory is simpler if $\varphi$ is a state. Then we may define the
embeddings $\kappa_p^{(\varphi)}: \HGP{p}{\varphi}{} \rightarrow
\HGP{1}{\varphi}{}: a \mapsto h^{1/(2p^*)}ah^{1/(2p^*)}$. The Kosaki
interpolation spaces (\cite{Kos2}) then correspond to the spaces $L^p(\M,
\varphi) = \kappa_p^{(\varphi)}(L^p_\varphi(\M))$ equipped with the norm
$\|a\|_p^\varphi = \|(\kappa_p^{(\varphi)})^{-1}a\|_p$. In this setting the
derivative $h$ may also be used to define embeddings $\M \rightarrow
L^p_\varphi(\M):a \mapsto h^{(1-c)/p}ah^{c/p}$ $(0 \leq c \leq 1)$ of $\M$ into
$L^p_\varphi(\M)$ (\cite{GL1}). For these embeddings the case $c = \frac{1}{2}$
has the added advantage of being positivity preserving, and so for this
distinguished case we will employ the notation $\mathfrak{i}^{(p)}$ for the
associated embedding.

As we have seen above, the Terp interpolation spaces are defined
only for the situation when $c=1/2$. The interested reader can find
a further generalization of the interpolation for the weight case, so
as to incorporate the cases when $c\neq 1/2$, in \cite{I}.  In
settings where several algebras or weights are involved we will
employ suitable subscripts to distinguish these cases.

In the sequel, by the term {\em Jordan ${}^*$-morphism} we understand a map
from a $C^*$-algebra into another $C^*$-algebra which preserves
adjoints and squares of elements.

\section{Defining generalised composition operators}

Let $(X_i, \Sigma_i, m_i)$ $(i = 1, 2)$ be standard Borel spaces
and let $T : Y \subset X_2 \rightarrow X_1$ be a given
non-singular measurable transformation from a measurable subset
$Y$ of $X_2$ into $X_1$. Then for $\infty > p \geq q \geq 1$ the
formula $C_T(f)(t) = f \circ T(t)$ if $t \in Y$ and $C_T(f)(t) =
0$ otherwise, induces a bounded linear operator $C_T: \LP{p}{1}
\rightarrow \LP{q}{2}$ if and only if $m_2 \circ T^{-1}$ is
absolutely continuous with respect to $m_1$ and
$\frac{\mathrm{d}m_2 \circ T^{-1}}{\mathrm{d}m_1}$ belongs to
$\LP{r}{1}$ where $r = \frac{p}{p-q}$ \cite{L1}. So we see that
when it comes to the formal existence of a (generalised)
composition operator in the case where $1 \leq q < \infty$, some
form of absolute continuity is crucial. (Boundedness of the
composition operator is in turn conditioned by the behaviour of
the associated Radon-Nikodym derivative.) We will see that even in
the noncommutative world it is precisely some form of absolute
continuity that once again enables us to formally introduce the
concept of a (generalised) composition operator.

Let $\M$ be a von Neumann algebra with $fns$ weight $\varphi$ and let $h =
\frac{\mathrm{d} \widetilde{\varphi}}{\mathrm{d}\tau}$. Then the span of the
set $$\{h^{1/(2p)}eh^{1/(2p)} | e \in \M \mbox{ a projection}, \varphi_1(e) <
\infty\}$$ is known to be norm dense in $L^p_{\varphi}(\M)$ if $1 \leq p <
\infty$. We may think of this span as representing the \emph{simple functions}
in $\HGP{p}{\varphi}{}$. In the context of classical $L^p$ spaces on standard
Borel measure spaces a bounded linear operator from $L^p$ to $L^q$ is known to
be a (generalised) composition operator precisely if it takes characteristic
functions in $L^p$ to characteristic functions in $L^q$. (See for example
\cite{L1}.) Now let $h_i = \frac{\mathrm{d}
\widetilde{\varphi_i}}{\mathrm{d}\tau_i}$, and let $J: \M_1 \rightarrow \M_2$
be a normal Jordan $*$-morphism satisfying the condition that for any
projection $e \in \M_1$ with $\varphi_1(e) < \infty$, we always have that
$\varphi_2(J(e)) < \infty$. In such a case the formal process
$h_1^{1/(2p)}ah_1^{1/(2p)} \rightarrow h_2^{1/(2q)}J(a)h_2^{1/(2q)}$ $(a \in
\M)$ is at least densely defined on $\HGP{p}{\varphi}{1}$. If indeed the
process extends to a bounded map $C_J : L^p_{\varphi_1}(\M_1) \rightarrow
L^q_{\varphi_2}(\M_2)$, then by analogy with the classical context mentioned
above, we may think of $C_J$ as a (generalised) composition operator induced by
$J$. We proceed to indicate that the condition regarding the Jordan
$*$-morphism's action on projections with finite weight may be interpreted as a type of local
absolute continuity. Thus the proposed definition of composition operators
compares well with the classical setting in that here too some form of absolute
continuity of $\varphi_2 \circ J$ with respect to $\varphi_1$ is a prerequisite
for the existence of a composition operator.

We start with a simple generalisation of a well known fact regarding
absolute continuity of finite measures. First, we give the following
definitions.

\begin{definition}
Let $\varphi_0,\varphi_1$ be weights on a von Neumann algebra $\M$.
\begin{enumerate}
\item
We say that $\varphi_0$ is $\epsilon$-$\delta$ absolutely continuous with
respect to $\varphi_1$ if, for every $\epsilon
> 0$ we can then find a $\delta > 0$ so that for any projection $e
\in \M$ with $\varphi_1(e) < \delta$ we will have that $\varphi_0(e) <
\epsilon$. For a projection $e\in M$, the weight $\varphi_0$ is called
$\epsilon$-$\delta$ absolutely continuous with respect to $\varphi_1$ on $e$ if
the restriction of $\varphi_0$ to the von Neumann algebra $eMe$ is
$\epsilon$-$\delta$ absolutely continuous with respect to the restriction of
$\varphi_1$ to $eMe$.
\item
We say that $\varphi_0$ is locally absolutely continuous with respect to
$\varphi_1$ if, for each projection $e\in M$, $\varphi_1(e)<\infty$ implies
$\varphi_0(e)<\infty$. If this is the case, we write $\varphi_0 \ll_{loc}
\varphi_1$.
\end{enumerate}
\end{definition}

We are going to show that (under very mild conditions) local
absolute continuity is, in fact, absolute continuity on each
projection of finite weight, so that the name is well chosen. In the
sequel, we assume that the weight $\varphi_1$ is semifinite.
Although this assumption is not really needed, it makes statements
of the results slightly easier, and is exactly what we need in
practice. Moreover, if $\varphi_1$ is not semifinite, there exists a
greatest projection $e$ such that $\varphi_1$ is semifinite when
restricted to $eMe$; it is enough to take for $e$ the unit of the
von Neumann algebra generated by projections of finite
$\varphi_1$-weight. Thus, we can always restrict our attention to an
algebra on which the weight in question is semifinite.

\begin{lemma} \label{ncfac}
Let $\varphi_1, \varphi_0$ be normal states on a von Neumann algebra $\M$ with
$\varphi_1$ also faithful. Then $\varphi_0$ is $\epsilon$-$\delta$ absolutely
continuous with respect to $\varphi_1$.
\end{lemma}

\begin{proof}
Note that the sets $\{x\in (\M)_1\colon \varphi_1(x^*x)<\epsilon\}$,
with $\epsilon>0$, form a basis of neighbourhoods of zero for the
strong topology on the unit ball of $\M$. Hence, the conclusion
follows from the strong continuity of $\varphi_0$ on the ball.
\end{proof}

The next two lemmas collect various facts belonging to the
mathematical folklore.

\begin{lemma}\label{nonat}
Let $\M$ be a von Neumann algebra with no minimal projections.
Then any maximal abelian von Neumann subalgebra $\M_0$ of $\M$
also has no minimal projections \cite{GJL}. If $\M$ admits of a faithful
normal state $\varphi$, then the algebra $\M_0$ corresponds to a
classical $L^\infty(\Omega, \Sigma, \mu_\varphi)$, where $(\Omega,
\Sigma, \mu_\varphi)$ is a nonatomic probability space and the
measure $\mu_\varphi$ is defined by $\mu_\varphi(E) =
\varphi(\chi_E)$ for each $E \in \Sigma$.
\end{lemma}

\begin{proof}
The first statement was noted in \cite{GJL}. The second follows 
from the fact that any commutative von Neumann subalgebra $\M_0$ will
correspond to some $L^\infty(\Omega, \Sigma, \nu)$. In particular given 
a faithful normal state $\varphi$ on $\M$, it is an exercise to
show that the restriction of $\varphi$ to $\M_0 = L^\infty(\Omega,
\Sigma, \nu)$ defines a probability measure $\mu_\varphi = \mu$ on $(\Omega,
\Sigma)$ (with the same sets of measure zero as $\nu$) by means of
the prescription $\mu(E) = \varphi(\chi_E)$ $E \in \Sigma$. Replacing 
$\nu$ by $\mu$ if necessary, all that remains is to note that the
subalgebra $\M_0 = L^\infty(\Omega, \Sigma, \mu)$ has no minimal
projections precisely when $(\Omega, \Sigma, \mu)$ is nonatomic.
\end{proof}

\begin{lemma} \label{nonat2}
Let $(\Omega, \Sigma, \mu)$ be a nonatomic probability space and let $\nu$ be a
measure on $(\Omega, \Sigma)$ which is $\epsilon$-$\delta$ absolutely
continuous with respect to $\mu$. Then $\nu$ is a finite measure.
\end{lemma}

\begin{proof}
Let $\epsilon$ be given and select $\delta$ so that for any $E \in
\Sigma$ with $\mu(E) < \delta$ we will have that $\nu(E) <
\epsilon$. We show that we may write $\Omega$ as the union of a
finite collection $E_1, E_2, \dots, E_n$ of disjoint sets in
$\Sigma$ with $\mu(E_k) < \delta$ for each $1 \leq k \leq n$. It
then trivially follows that $\nu(\Omega) = \sum_{k=1}^n \nu(E_k) <
n \epsilon < \infty$ as required. To see that such a partitioning
of $\Omega$ is indeed possible let $n \in \mathbb{N}$ be given
such that $\frac{1}{n} < \delta$, and use Zorn's lemma to find a
maximal set $E_1 \in \Sigma$ with $\mu(E_1) \leq \frac{1}{n}$. Now
given any $E \in \Sigma$ with $\mu(E) < \frac{1}{n}$ we can then
use the nonatomicity of $(\Omega, \Sigma, \mu)$ to find a larger
set $F \in \Sigma$ with $\mu(E) < \mu(F) < \frac{1}{n}$. Hence the
maximality of $E_1$ ensures that $\mu(E_1) = \frac{1}{n}$. To
complete the proof we may now continue inductively by finding a
measurable subset $E_2$ of $\Omega - E_1$ such that $\mu(E_2) =
\frac{1}{n-1}\mu(\Omega - E_1) = \frac{1}{n}$, and so on.
\end{proof}

\begin{theorem}\label{lac}
Let $\M$ be an arbitrary von Neumann algebra, $\varphi_1$ be  a
faithful normal semifinite weight on $\M$ and $\varphi_0$ a normal
weight on $\M$, semifinite on its atomic part. The following
conditions are equivalent:
\begin{enumerate}
\item $\varphi_0$ is locally absolutely continuous with respect to
$\varphi_1$;
\item $\varphi_0$ is $\epsilon$-$\delta$ absolutely continuous with respect to
$\varphi_1$ on each projection $e\in M$ with $\varphi_1(e)<\infty$.
\end{enumerate}
\end{theorem}

\begin{proof} That local absolute continuity implies $\epsilon$-$\delta$ absolute continuity
on each projection $e\in M$ with $\varphi_1(e)<\infty$ follows
immediately from Lemma \ref{ncfac}. For the reverse implication, we
fix $\epsilon>0$ and take the corresponding $\delta$ from the
definition of the $\epsilon$-$\delta$ absolute continuity. One notes
first that if the algebra $\M$ is a direct sum of a finite number of
summands, it is enough to prove the implication on each summand
separately. Thus, it is enough to consider the following four cases:
\begin{enumerate}
\item The algebra $\M$ is non-atomic (i.e. it has no minimal
projections). Let $e$ be a projection such that
$\varphi_1(e)<\infty$. Now since $e$ belongs to a maximal abelian
subalgebra, say $\M_0$, of $e\M e$, it suffices to prove that if
$\varphi_0$ restricts to a normal weight on $\M_0$ which is
$\epsilon-\delta$ absolutely continuous with respect to the action
of $\varphi_1$ on $\M_0$, then $\varphi_0|_{\M_0}$ is a finite
weight on $\M_0$. Without loss of generality we may of course
normalise the action of $\varphi_1$ on $\M_0$. Then by Lemma
\ref{nonat} the algebra $\M_0$ corresponds to a classical
$L^\infty(\Omega, \Sigma, \mu)$, where $(\Omega, \Sigma, \mu)$ is a
nonatomic probability space and the measure $\mu$ is defined by
$\mu(E) = \varphi_1(\chi_E)$ for each $E \in \Sigma$. In a similar
fashion the weight $\varphi_0$ also defines a measure $\nu$ on
$(\Omega, \Sigma)$ by means of the formula $\nu(E) =
\varphi_0(\chi_E)$ for each $E \in \Sigma$. We may then directly
conclude from Lemma \ref{nonat2} that $\varphi_0(e) =
\varphi_0(\chi_{\Omega}) = \nu(\Omega) < \infty$ as required.
\item The algebra $\M$ is a factor of type $I_\infty$ (where $\infty$
stands for any infinite cardinal). Let $e$ be a projection such that
$\varphi_1(e)=1$. Note that $\varphi_0$ is finite on any minimal
projection of $\M$, by semifiniteness. Hence we may assume that $e$
is (properly) infinite. Write $e$ in the form $e=\Sigma_{k=1}^\infty
e_k$, where the projections $e_k$ are all equivalent to $e$. Choose
$n$ so that $\frac{1}{n}<\delta$. Then, for some $k$,
$\varphi_1(e_k)<1/n$. Since $\varphi_1(e_k)\leq\varphi_1(e-e_k)$ and
$e_k\sim e-e_k$, there is a projection $f_1$ in $\M$ such that
$f_1\leq e$ and $\varphi_1(f_1)=1/n$ (see, for example, \cite{GP},
Proposition 1.1). The rest follows the lines of the proof of Lemma
\ref{nonat2}.
\item Assume now that $\M$ is finite and atomic. Then $\M$ is of the form
$\Sigma_{i\in I}^\oplus M_i$, where each $M_i$ is a factor of type
$I_{n_i}$ with $n_i<\infty$. As before, since $\varphi_0$ is
semifinite on $\M$, we may assume that the index set $I$ is
infinite. Let $e$ be a projection such that $\varphi_1(e)<\infty$.
Then $e$ is of the form $\sum_{i\in I}e_i$ and there exists a finite
subset $J$ of $I$ such that $\varphi_1(\sum_{i\in I\setminus
J}e_i)<\delta$. Hence $\varphi_0(e)=\sum_{i\in
J}\varphi_0(e_i)+\varphi_0(\sum_{i\in I\setminus J}e_i) < \infty$,
by the $\epsilon$-$\delta$ condition and the semifiniteness of
$\varphi_0$.
\item Assume finally that $\M$ is an infinite direct sum of type
$I_\infty$ factors. We obtain the result as in (3), from the
$\epsilon$-$\delta$ condition and (2).
\end{enumerate}

Note that we did not use the assumption that $\varphi_0$ is normal
in the proof of the reverse implication.
\end{proof}

\begin{remark}\label{semif} Let $\M$ be a von Neumann algebra with two normal
weights $\varphi_0$ and $\varphi_1$, with $\varphi_1$ also
semifinite and faithful. Now if $\varphi_0$ was locally absolutely
continuous with respect to $\varphi_1$, then $\varphi_0$ would
also be semifinite! To see this all we need to notice is that the
linear span of all projections $e\in \M$ with $\varphi_1(e) <
\infty$ is $\sigma$-weakly dense in $\M$.
\end{remark}

In the sequel, whenever we deal with a von Neumann algebra $\M$ with a fixed
weight $\varphi$, we shall write $\M^{(0)}$ for the span of the set $\{e | e
\in \M \mbox{ a projection}, \varphi(e) < \infty\}$. The
weight used to define $\M^{(0)}$ will always be clear from the context.

We are now ready to formally define the concept of a composition
operator on Haagerup $L^p$-spaces.

\begin{definition}
Let $J: \M_1 \rightarrow \M_2$ be a normal Jordan $*$-morphism, let $h_i =
\frac{\mathrm{d} \widetilde{\varphi_i}}{\mathrm{d}\tau_i}$. Given $1 \leq q
\leq p < \infty$, we say that $J$ induces a \emph{generalised composition
operator} (or just a \emph{composition operator} if $J(\I) = \I$) from
$L^p_{\varphi_1}(\M_1)$ into $L^q_{\varphi_2}(\M_2)$
\begin{itemize}
\item if $\varphi_2 \circ J$ is locally absolutely
continuous with respect to $\varphi_1$,
\item and if the process $h_1^{1/(2p)}ah_1^{1/(2p)} \rightarrow
h_2^{1/(2q)}J(a)h_2^{1/(2q)}$ $(a \in \M_1^{(0)})$ is continuous.
\end{itemize} The above process then extends uniquely to a bounded
map $C_J : L^p_{\varphi_1}(\M_1) \rightarrow L^q_{\varphi_2}(\M_2)$,
which we shall call the \emph{(generalised) composition operator}
induced by $J$ from $L^p_{\varphi_1}(\M_1)$ into
$L^q_{\varphi_2}(\M_2)$. (Here we used the fact that
$\mathfrak{i}^{(p)}(\M_1^{(0)})$ is norm dense in
$L^p_{\varphi_1}(\M_1)$.)
\end{definition}

\begin{remark}
By analogy with the above definition we may say that $J$ induces a generalised
composition operator from $\M_1 = L^\infty_{\varphi_1}(\M_1)$ into
$L^q_{\varphi_2}(\M_2)$ ($1 \leq q < \infty$) if the map $C_J: a \rightarrow
h_2^{1/(2q)}J(a)h_2^{1/(2q)}$ $(a \in \M)$ is well-defined and continuous from
$\M_1$ into $L^q_{\varphi_2}(\M_2)$. For this map to be well-defined we at
least need $J(\I)h_2^{1/(2q)}$ to be closable with closure an element of
$L^{2q}_{\varphi_2}(\M_2)$ (see (\cite{GL2}; 2.7 and 2.8).

Conversely if $J(\I)h_2^{1/(2q)}$ is closable with closure an element of
$L^{2q}_{\varphi_2}(\M_2)$, then the above map is well-defined and continuous.
In fact, for for any $a \in \M_1$ we will have that $J(a)J(\I) \in
\mathfrak{n}_{\varphi_2}^{(2q)}$ and hence that
$[J(\I)h_2^{1/(2q)}]^*[J(a)J(\I)h_2^{1/(2q)}] =
[h_2^{1/(2q)}J(\I)J(a)J(\I)h_2^{1/(2q)}] = [h_2^{1/(2q)}J(a)h_2^{1/(2q)}] \in
L^q_{\varphi_2}(\M_2)$. Now for any $a = a^* \in \M_1$ we have that
$-\|a\|_\infty\I \leq a \leq \|a\|_\infty\I$, and hence that
$$-\|a\|_\infty[h_2^{1/(2q)}J(\I)h_2^{1/(2q)}] \leq
[h_2^{1/(2q)}J(a)h_2^{1/(2q)}] \leq
\|a\|_\infty[h_2^{1/(2q)}J(\I)h_2^{1/(2q)}].$$ From this it follows that
$$\|[h_2^{1/(2q)}J(a)h_2^{1/(2q)}]\|_q \leq \|a\|_\infty
\|[h_2^{1/(2q)}J(\I)h_2^{1/(2q)}]\|_q$$ for each $a = a^* \in \M_1$ (or rather
$\|C_J(a)\|_q \leq \|C_J(\I)\|_q \|a\|_\infty)$. This clearly suffices to force
continuity of the induced map.
\end{remark}

\section{Identifying and describing composition operators}

Having introduced the concept of a composition operator on
Haagerup $L^p$-spaces we now focus on the two-fold task of firstly
finding a way to identify those operators that actually are
composition operators, and secondly describing those Jordan
morphisms between von Neumann algebras that do indeed induce
composition operators on the associated $L^p$-spaces.

\subsection*{Operators on Haagerup $L^p$-spaces that come from
Jordan $*$-morphisms}

We noted earlier that on classical $L^p$ spaces of standard Borel
measure spaces, bounded linear operators from $L^p$ to $L^q$ are
(generalised) composition operators precisely when they take
characteristic functions in $L^p$ to characteristic functions in
$L^q$. (See for example \cite{L1}.) The primary result of this
section shows that a similar structure pertains even in the
noncommutative context. In this regard we remind the reader that
given a von Neumann algebra $\M$ equipped with a faithful normal
semifinite weight $\varphi$, the role that is classically played by
characteristic functions in $L^p$ will here be played by elements
of the form $h^{1/(2p)}eh^{1/(2p)}$ where $e$ is a self-adjoint
projection in $\M$ with finite weight, and $h
=\frac{\mathrm{d}\widetilde{\varphi}}{\mathrm{d}\tau}$. Thus by
analogy with the classical setting, it is natural to try and
describe composition operators in terms of their action on
elements of the above form.

\begin{definition} Let $1 \leq p, q \leq \infty$ and let $\M_i$ $(i = 1, 2)$ be
von Neumann algebras equipped with faithful normal semifinite weights
$\varphi_i$. We say that a bounded linear operator $S:L^p_{\varphi_1}(\M_1)
\rightarrow L^q_{\varphi_2}(\M_2)$ \emph{preserves characteristic functions} if
for any projection $e \in \M_1$ (with $\varphi_1(e) < \infty$ if $p <\infty$)
there exists a unique projection $\widetilde{e} \in \M_2$ (with
$\varphi_2(\widetilde{e}) < \infty$ if $q < \infty$) such that
$S(h_1^{1/(2p)}eh_1^{1/(2p)}) = h_2^{1/(2q)}\widetilde{e}h_2^{1/(2q)}$ (where
$h_i =\frac{\mathrm{d}\widetilde{\varphi_i}}{\mathrm{d}{\tau_i}})$.
\end{definition}

\begin{theorem} Let $1 \leq p, q \leq \infty$ and let $\M_i$ $(i = 1, 2)$ be von
Neumann algebras equipped with faithful normal semifinite weights
$\varphi_i$. Let ${\mathcal C}(\M_1)$ denote the $C^*$-subalgebra
of $\M_1$ generated by $\M_1^{(0)}$. Let $S:L^p_{\varphi_1}(\M_1)
\rightarrow L^q_{\varphi_2}(\M_2)$ be a bounded linear operator
which preserves characteristic functions.

If $p = \infty$ then for some Jordan morphism $J: \M_1 \rightarrow \M_2$, $S$
is precisely of the form $a \rightarrow h_2^{1/(2q)}J(a)h_2^{1/(2q)}$ where $a
\in \M_1$.

If $p < \infty$, there exists a (not necessarily normal) Jordan $*$-morphism
$J: \mathcal{C}(\M_1) \to \M_2$ such that $S$ appears as the continuous
extension of the map $$h_1^{1/(2p)}ah_1^{1/(2p)} \rightarrow
h_2^{1/(2q)}J(a)h_2^{1/(2q)}$$where $a \in \M_1^{(0)}$. In this case $J$ will
be normal precisely when it satisfies the requirement that if mutually
orthogonal projections $e_1, \ldots \, e_n$ in $\M_1^{(0)}$, sets of mutually
orthogonal projections $\{f_i^{(k)}\}_I~ (k = 0, 1, \ldots, m)$ in
$\M_1^{(0)}$, and positive scalars $\lambda_1, \lambda_2, \ldots, \lambda_n$
and $\mu_1, \mu_2, \ldots, \mu_m$ are such that $$\sum^n_{p=1} \lambda_{p}e_p
\leq \sum^m_{k=1} \mu_{k}\left(\sum_{i \in I} f_i^{(k)} \right),$$it
then follows that $$\sum^n_{p=1} \lambda_p J(e_p) \leq \sum^m_{k=1}
\mu_k \left( \sum_{i \in I} J(f_i^{(k)}) \right).$$

If in fact $\varphi_1, \varphi_2$ are states and $p < \infty$,
then $J$ is necessarily normal (and of course defined on all of
$\M_1$).
\end{theorem}

Note that in the above we do not require that $q \leq p$.

Next let $1 \leq p < \infty$. In the commutative case normality of
$J$ will then still be automatic even if $\varphi_1, \varphi_2$
are not states. (See \cite[4.3 \& 4.15(iii)]{L1}.) However,
although we have no proof for this as yet, we suspect that in the
noncommutative setting $\sigma$-finiteness is essential to obtain
automatic normality of $J$.

\begin{proof}
The proofs for the cases $p = \infty$ and $p < \infty$ are similar,
and hence we prove only the latter case. Let
$S:L^p_{\varphi_1}(\M_1) \rightarrow L^q_{\varphi_2}(\M_2)$ be a
bounded linear operator which preserves characteristic functions in
the sense described above.

Firstly note that by hypothesis $S$ will map all elements of the form
$h_1^{1/(2p)}ah_1^{1/(2p)}$ where $a \in \M_1^{(0)}$ onto elements of the form
$h_2^{1/(2q)}\widetilde{a}h_2^{1/(2q)}$ where $\widetilde{a} \in \M_2^{(0)}$.
So if for some $a \in \M_1^{(0)}$ we have $S(h_1^{1/(2p)}ah_1^{1/(2p)}) =
h_2^{1/(2q)}\widetilde{a}h_2^{1/(2q)}$, we set $J(a) = \widetilde{a}$. The
linearity of $S$ and the injectivity of
$\mathfrak{i}^{q}$ ensures that $J: \M_1^{(0)} \to \M_2$ is well-defined and linear.

Notice that if $e$ and $f$ are mutually orthogonal projections in
$\M_1^{(0)}$, then by construction each of $J(e)$, $J(f)$ and $J(e +
f) = J(e) + J(f)$ is also a projection. However the latter can only
hold if in fact $J(e) \perp J(f)$. It therefore follows that $J$
preserves the orthogonality of projections in $\M_1^{(0)}$. But then
$J$ will also preserve the order of projections.

Now let $a \in \M_1^{(0)}$ be given with $a = a^*$. Since $a$ is in
$\M_1^{(0)}$, we surely have $\varphi_1(\supp(a)) < \infty$. For the
sake of simplicity write $e = \supp(a)$. Then by passing to Riemann
sums of spectral projections of $a$, we can find a sequence
$$b_n = \sum_{k=1}^{m_n} \mu_k^{(n)}e_k^{(n)} \in \M_1^{(0)}$$
converging uniformly to $a$ such that for each fixed $n \in
\mathbb{N}$:
\begin{itemize}
\item the projections $\{e_k^{(n)} | 1 \leq k \leq m_n\}$ are
mutually orthogonal and satisfy $0 \leq e_k^{(n)} \leq e$;
\item $-\|a\| \leq \mu_k^{(n)} \leq \|a\|$,  $1 \leq k \leq m_n$.
Then of course $-\|a\|e \leq b_n \leq \|a\|e$.
\end{itemize}
Since $J$ preserves both the order and orthogonality of
projections, it is clear from the above facts that $$-\|a\|J(e)
\leq J(b_n) \leq \|a\|J(e)$$ and hence that
$$-\|a\|h_2^{1/(2q)}J(e)h_2^{1/(2q)} \leq h_2^{1/(2q)}J(b_n)h_2^{1/(2q)}
\leq \|a\|h_2^{1/(2q)}J(e)h_2^{1/(2q)}$$ for each $n$. Since $eh_1^{1/(2p)}$ is
measurable and $\supp(b_n)\leq e$, the uniform convergence of the $b_n$'s to
$a$ ensures that $h_1^{1/(2p)}b_nh_1^{1/(2p)} \to h_1^{1/(2p)}ah_1^{1/(2p)}$.
The continuity of $S$ then yields $$h_2^{1/(2q)}J(b_n)h_2^{1/(2q)} =
S(h_1^{1/(2p)}b_nh_1^{1/(2p)}) \to S(h_1^{1/(2p)}ah_1^{1/(2p)}) =
h_2^{1/(2q)}J(a)h_2^{1/(2q)}$$ in $L^q(\M_2)$. Together these two facts force
$$-\|a\|h_2^{1/(2q)}J(e)h_2^{1/(2q)} \leq h_2^{1/(2q)}J(a)h_2^{1/(2q)} \leq
\|a\|h_2^{1/(2q)}J(e)h_2^{1/(2q)}.$$
Let $b \in \mathfrak{m}_{\varphi_2, +}$ be given. Since $\mathfrak{i}^{(q^*)}(b)
\in L^{q^*}_+(\M_2)$, we have
$$-\|a\|\tr(\mathfrak{i}^{(q^*)}(b)\mathfrak{i}^{(q)}(J(e))) \leq
tr(\mathfrak{i}^{(q^*)}(b)\mathfrak{i}^{(q)}(J(a))) \leq
\|a\|\tr(\mathfrak{i}^{(q^*)}(b)\mathfrak{i}^{(q)}(J(e)))$$or equivalently
$$-\|a\|\tr(\mathfrak{i}^{(1)}(b)J(e)) \leq tr(\mathfrak{i}^{(1)}(b)J(a)) \leq
\|a\|\tr(\mathfrak{i}^{(1)}(b)J(e)).$$On applying \cite[Proposition 2.11(b)]{GL2},
it now follows that
$-\|a\|J(e) \leq J(a) \leq \|a\|J(e)$, and hence that $\|J(a)\| \leq \|a\|$. Thus $J$ is
bounded. By continuity we may then extend $J$ to the uniform closure of
$\M_1^{(0)}$. This closure is however exactly $\mathcal{C}(\M_1)$. To see this
note that if $b = b^*$ is in the dense $*$-subalgebra of $\mathcal{C}(\M_1)$
generated by finite algebraic combinations of elements of $\M_1^{(0)}$, then
$\varphi_1(\supp(b)) < \infty$, and hence as before by passing to Riemann sums
we may write $b$ as a norm limit of terms of the form $d_n = \sum_{k=1}^{m_n}
\mu_k^{(n)}f_k^{(n)} \in \M_1^{(0)}$ where the $f_k^{(n)}$'s are mutually
orthogonal. Then $b = \lim_n d_n \in \overline{M_1^{(0)}}$. Thus
$\mathcal{C}(\M_1) \subset \overline{\M_1^{(0)}}$. The converse inclusion is
clear. Now with $b$ as above, notice that also $J(b^2) = \lim_n J(d_n^2) =
\lim_n J(\sum_{k=1}^{m_n} (\mu_k^{(n)})^2 f_k^{(n)}) = \lim_n(\sum_{k=1}^{m_n}
\mu_k^{(n)}J(f_k^{(n)}))^2 = \lim_n(J(d_n)^2) = J(b)^2$. Thus $J$ preserves
squares of self-adjoint elements on $\mathcal{C}(\M_1)$, and hence must be a
Jordan $*$-morphism.

The claim about the normal extension of $J$ to all of $\M_1$ may
be proved by a similar argument as was employed in the proof of
the implications $(iii) \Rightarrow (iv) \Rightarrow (ii)$
in \cite[4.4]{L1}. The only change that needs to be made is that
wherever semifiniteness of $\M_1$ was used in \cite{L1} to select
a finite subprojection $e$, we should here use the semifiniteness
of $\varphi_1$ to select a subprojection $e$ with $\varphi_1(e) <
\infty$.

It remains to show that $J$ is normal when $\varphi_1, \varphi_2$ are states
and $p < \infty$. Since $\varphi_1$  is a state, it is clear that in this case
$J$ is defined on all of $\M_1$. So suppose that $p < \infty$, and let
$\{e_\mu\}_{\mu}$ be a set of mutually orthogonal projections in $\M_1$. If we
can show that $J(\sum_{\mu} e_{\mu}) = \sum_{\mu} J(e_{\mu})$, $J$ will be
normal by \cite[4.3]{L1}. Now $$e = \sum_{\mu} e_{\mu}$$ is of course a
projection in $\M_1$ with convergence of the series taking place in the
$\sigma$-strong topology (and hence also the weak* topology) of $\M_1$. But
then $$h_1^{1/(2p)}eh_1^{1/(2p)} = \sum_{\mu} h^{1/(2p)}e_{\mu}h^{1/(2p)}$$
with convergence taking place in the weak topology of $L^p(\M_1)$. To see this
note that if $a_\lambda \to a$ in the weak* topology on $\M_1$, then for any $b
\in L^{p^*}_{\varphi_1}(\M_1)$ we will have $\tr((h^{1/(2p)}a_\lambda
h^{1/(2p)})b) = \tr(a_\lambda(h^{1/(2p)}b h^{1/(2p)})) \to \tr(a(h^{1/(2p)}b
h^{1/(2p)})) = \tr((h^{1/(2p)}a h^{1/(2p)})b)$ (since then
$h^{1/(2p)}bh^{1/(2p)} \in L^{1}_{\varphi_1}(\M_1)$).

Since $S$ is norm continuous, it is also weak-weak continuous.
Thus we may conclude from the above that
$$h_2^{1/(2q)}J(e)h_2^{1/(2q)} = \sum_{\mu}
h_2^{1/(2q)}J(e_{\mu})h_2^{1/(2q)}$$ with convergence taking place in the weak
topology on $L^q_{\varphi_2}(\M_2)$. But since $\{J(e_\mu)\}_{\mu}$ is a set of
mutually orthogonal projections in $\M_2$, it follows that $$f = \sum_{\mu}
J(e_{\mu})$$ is a projection in $\M_2$ with convergence taking place in the
weak* topology on $\M_2$. Now if $q = \infty$, uniqueness of limits will then
force $J(\sum_{\mu} e_{\mu}) = J(e) = f = \sum_{\mu} J(e_{\mu})$. If however $q
< \infty$, we may argue as before to conclude that $$h_2^{1/(2q)}fh_2^{1/(2q)}
= \sum_{\mu} h_2^{1/(2q)}J(e_{\mu})h_2^{1/(2q)}$$ with convergence taking place
in the weak topology on $L^q_{\varphi_2}(\M_2)$. Once again uniqueness of
limits will then force $h_2^{1/(2q)}J(e)h_2^{1/(2q)} =
h_2^{1/(2q)}fh_2^{1/(2q)}$. Since
$h_2$ is an injective positive element of
$L^1_{\varphi_2}(\M_2)$, this is enough to ensure that $J(\sum_{\mu} e_{\mu}) =
J(e) = f = \sum_{\mu} J(e_{\mu})$ as required.
\end{proof}

\subsection*{Jordan $*$-morphisms that induce operators on
Haagerup $L^p$-spaces}

The main focus of this subsection is to try and describe those Jordan
$*$-morphisms which allow for the construction of a (generalised) composition
operator on a given pair of $L^p$-spaces. Although we do not succeed in giving
a completely general description, we do manage to describe a large class of
morphisms from which we may construct such operators. We will assume throughout
that $\M_i$ ($i = 1, 2$) are von Neumann algebras equipped with faithful normal
semifinite weights $\varphi_i$, and that $J: \M_1 \to \M_2$ is a normal Jordan
$*$-morphism. Moreover, $\mathcal{B}$ is the von Neumann algebra generated by
$J(\M_1)$ and $\varphi_{\mathcal{B}}$ denotes the restriction of $\varphi_2$ to
$\mathcal{B}$. Note that the unit of $\mathcal{B}$ is $J(\I)$.

It turns out that the construction of \emph{composition operators} from such a
Jordan $*$-morphism may be broken up into five distinct steps. To avoid any
pathologies associated with this process, we will for the remainder of this
section consistently assume that $\varphi_2 \circ J$ is locally absolutely
continuous with respect to $\varphi_1$. To gain some clarity regarding the
processes involved, we first take some time to review the classical situation.

\subsubsection*{Preamble to the construction of
composition operators}

Let $(X_i, \Sigma_i, m_i)$ $(i = 1, 2)$ be measure spaces and let
$T : Y \subset X_2 \rightarrow X_1$ be a given non-singular
measurable transformation from a measurable subset $Y$ of $X_2$
into $X_1$. For any $q$ we may then regard $L^q(Y, m_2)$ as a
subspace of $L^q(X_2, m_2)$ by simply assigning the value $0$ on
$X_2 \backslash Y$ to each element of $L^q(Y, m_2)$. If the
process $f \rightarrow f \circ T$ \emph{directly} yields a bounded
linear operator from $L^p(X_1, m_1)$ to $L^q(Y, m_2) \subset
L^q(X_2, m_2)$, we call the resultant operator a \emph{generalised
composition operator from $L^p(X_1, m_1)$ to $L^q(X_2, m_2)$} and
denote it by $C_T$ . If in fact $Y = X_2$, we simply call $C_T$ a
composition operator.

Notice that we may use $T$ to define a new measure $m_2 \circ
T^{-1}$ on $X_1$. With this new measure in place one should now be
very careful about what one calls a ``composition operator''. For
example the map $L^q(X_1, m_2 \circ T^{-1}) \rightarrow L^q(X_2,
m_2)$ defined by $f \rightarrow f \circ T$ is a very nice map (in
fact an isometry), but it is not a composition operator \emph{from
$L^p(X_1, m_1)$ to $L^q(X_2, m_2)$} in the true sense of the word.
Part of the problem is that the measure on the domain space is
wrong.

Now if we do have a bounded map of the form $C_T : L^p(X_1, m_1)
\rightarrow L^q(Y, m_2) \subset L^q(X_2, m_2): f \mapsto f\circ T$,
the construction of such a map may be broken up into five
subprocesses. In the following let $Z \in \Sigma_1$ be the support
of $m_2 \circ T^{-1}$ in $X_1$, let $\Sigma_2^Y = \{E \in \Sigma_2 |
E \subset Y\}$, and let $\Sigma_T$ be the $\sigma$-subalgebra of
$\Sigma_Y$ generated by sets of the form $T^{-1}(E)$ where $E \in
\Sigma_1$. Our composition operator is then made up of the following
processes:
\begin{enumerate}
\item[(I)] Restricting to the support of $m_2 \circ T^{-1}$:
$L^p(X_1, m_1) \rightarrow L^p(Z, m_1|_Z): f \mapsto f|_Z$

\item[(II)] Changing weights: $L^p(Z, m_1|_Z) \rightarrow L^q(Z, m_2 \circ
T^{-1}): f \mapsto f$

\item[(III)] Isometric equivalence of spaces: $L^q(Z, \Sigma_1^Z, m_2 \circ T^{-1})
\rightarrow L^q(Y, \Sigma_T, m_2): f \mapsto f \circ T$ (Here
$\Sigma_1^Z = \{E \in \Sigma_1 | E \subset Z\}$.)

\item[(IV)] Refining the $\sigma$-algebra: $L^q(Y, \Sigma_T, m_2) \rightarrow
L^q(Y, \Sigma_2^Y, m_2): f \mapsto f$

\item[(V)] Canonical embedding: $L^q(Y, \Sigma_2^Y, m_2) \rightarrow
L^q(X_2, \Sigma_2, m_2): f \mapsto j(f)$ where $j(f) = f$ on $Y$
and $j(f) = 0$ on $X_2\setminus Y$.
\end{enumerate}

Notice that the map in step (V) will be the identity whenever
$X_2\setminus Y$ is a set of measure zero. Now for the combination
of these five processes to yield a composition operator, we must
careful about HOW we change weights. Suppose by way of example that
$m_1$ and $m_2 \circ T^{-1}$ have the same sets of measure zero and
that $\frac{\mathrm{d} m_1}{\mathrm{d}m_2 \circ T^{-1}}$ exists.
Then the map $f \rightarrow f(\frac{\mathrm{d} m_1}{\mathrm{d}m_2
\circ T^{-1}})^{1/p}$ will certainly yield an isometry from
$L^p(X_1, m_1)$ to $L^p(X_1, m_2 \circ T^{-1})$, but using this to
change weights will not in general yield a composition operator. In
the following we give some indication of how one may construct
``composition operators'' on noncommutative $L^p$-spaces associated
with von Neumann algebras, by successively extending each of these
processes to the noncommutative context. Thus given von Neumann
algebras $\M_i$ ($i = 1, 2$) the basic idea is to classify and study
those Jordan $*$-morphisms $J:\M_1 \rightarrow \M_2$ that
canonically induce bounded linear operators $L^p(\M_1) \rightarrow
L^q(\M_2)$ along the lines suggested above. We proceed to look at
noncommutative versions of each of the above steps.

\subsubsection*{Step (I): Reducing matters to the case where
$J$ is injective}

Notice that $\varphi_2 \circ J$ defines a semifinite normal weight
on $\M_1$. So the noncommutative analogue of the first step would be
to pass from $(\M_1, \varphi_1)$ to $(e\M_1e, e \varphi_1 e)$, where
$e$ is the support projection of $\varphi_2 \circ J$, in a way that
allows us to compare the associated $L^p$-spaces. The object of this
exercise is basically to reduce matters to the case where $\varphi_2
\circ J$ is also faithful. We point out that no real information is
lost in making such a reduction since it follows from $J(\I) = J(e)$
that $J(a) = J(\I a\I) = J(\I)J(a)J(\I) = J(e)J(a)J(e) = J(eae)$ for
each $a \in \M_1$. It turns out that such a reduction is always
possible. We start with two easy lemmas concerning facts generally
known, which we chose to prove here for completeness.

Note that the algebra generated by $J(e\M e)$ is the same as the algebra
generated by $J(\M)$, that is $\mathcal{B}$.

Assume now that we have a von Neumann algebra $\M$ acting in a Hilbert space
$H$, with a \emph{fns} weight $\varphi$. If $e$ is a projection from $\M$, we
denote by $\varphi_e$ the restriction of  $\varphi$ to $e\M e$. Furthermore, we
denote by $\tau_e$ the canonical trace on the crossed product $(e\M
e)\rtimes_{\sigma^{\varphi_e}}\mathbb{R}$. Finally, we put
$\tilde{e}:=\pi_\varphi(e)$.

\begin{lemma}\label{fixed}
If the projection $e$ belongs to the subalgebra $\M_\varphi$ of fixed points
for the modular group of $\M$ with respect to an \emph{fns} weight $\varphi$
(in particular, when $e$ is central), then $L^p_{\varphi_e}(e\M e)$ consists of
operators from $\tilde{e}L^p_\varphi(\M)\tilde{e}$ restricted to
$L^2(\mathbb{R},eH)$. Moreover,
$\frac{\mathrm{d}\widetilde{\varphi}}{\mathrm{d}\tau}$ commutes with
$\tilde{e}$ and $\frac{\mathrm{d}\widetilde{\varphi_e}}{\mathrm{d}\tau_e}$ may
be identified with the restriction of
$\frac{\mathrm{d}\widetilde{\varphi}}{\mathrm{d}\tau}\tilde{e}$ to
$L^2(\mathbb{R},eH)$.
\end{lemma}

\begin{proof}
It is clear that the weight $\varphi_e$ is faithful, normal and semifinite, and
the modular group for the pair $(e\M e, \varphi_e)$ is the restriction to $e\M
e$ of the modular group for $(\M,\varphi)$. Similarly, one checks easily that
$\tilde{e}$ projects $L^2(\mathbb{R},H)$ onto $L^2(\mathbb{R},eH)$.
Consequently, the operators $\pi_{\varphi_e}(eae)$ with $a\in\M$ are just
$\tilde{e}\pi_\varphi(a)\tilde{e}$ restricted to $L^2(\mathbb{R},eH)$.
Similarly, $\lambda_{\varphi_e}(s)$ is just $\lambda_{\varphi}(s)$ restricted
to $L^2(\mathbb{R},eH)$. Hence $(e\M
e)\rtimes_{\sigma^{\varphi_e}}\mathbb{R}=\tilde{e}(\cM{}{\varphi})\tilde{e}$,
where the von Neumann algebra on the right hand side of the equation acts on
$L^2(\mathbb{R},eH)$. Now, if $(\theta_s)$ is the dual action on
$\cM{}{\varphi}$, then it restricts to the dual action on $(e\M
e)\rtimes_{\varphi_e}\mathbb{R}$, and
$\theta_s(\tilde{e}x\tilde{e})=\exp(-s/p)\tilde{e}x\tilde{e}$ for each $x\in
L^p_\varphi(\M)$, which implies the required equality. The final claim follows
from noting that the shift operators $\lambda_\varphi(s)$ commute with
$\tilde{e}$, and that that
$\frac{\mathrm{d}\widetilde{\varphi}}{\mathrm{d}\tau}$ (resp.
$\frac{\mathrm{d}\widetilde{\varphi_e}}{\mathrm{d}\tau_e}$) is the (positive)
generator of the unitary group $\lambda_{\varphi}(s), s\in\mathbb{R}$
(resp.$\lambda_{\varphi_e}(s), s\in\mathbb{R}$).
\end{proof}

\begin{remark}
The lemma shows that there is a \emph{natural} embedding of
$L^p_{\varphi_e}(e\M e)$ into $L^p_\varphi(\M)$, namely $x\mapsto
\tilde{e}x\tilde{e}$, and that the image of $L^p_{\varphi_e}(e\M e)$ under the
embedding is exactly $\tilde{e}L^p(\M)\tilde{e}$. In the sequel we stick to the
usual convention of identifying $e$ with $\tilde{e}$ and $L^p_{\varphi_e}(e\M
e)$ with $eL^p_\varphi(\M)e$.
\end{remark}

\begin{lemma}\label{supp}
The support of $\varphi_2\circ J$ is central.
\end{lemma}

\begin{proof}
Let $z$ be a central projection in $\mathcal{B}$ such that $a\mapsto zJ(a)$ is
a *-homomorphism and $a\mapsto (J(\I)-z)J(a)$ is a *-antihomomorphism. If
$J(a)=0$ for some $a\in \M$, then $J(ab)=zJ(a)J(b)+(J(\I)-z)J(b)J(a)=0$ and
similarly $J(ba)=0$. Hence the kernel of $J$ is a two-sided ideal,
$\sigma$-weakly closed because of $J$'s normality. Thus there exists a central
projection $e$ such that $\ker(J)=e\M$ (see \cite[Proposition II.3.12]{Tak}).
Now, it follows from Remark \ref{semif} that $\varphi_2\circ J$ is semifinite,
which shows that its support must be equal to $\I-e$.
\end{proof}

The above results show that the reduction to the support of
$\varphi_2\circ J$ is, in fact, multiplication by a central
projection. Since for any pair $(\M, \varphi)$, central
projections are automatically fixed points of the modular group of
$\M$ induced by $\varphi$ (in fact they are even central in
$\cM{}{\varphi}$), in the light of Lemma \ref{fixed} this
reduction is particularly simple.

\subsubsection*{Step (II): Changing weights}\label{chwt}

Let $J$ be as before and let $e$ be the support projection of $\varphi_2 \circ
J = \varphi_J$. Our primary interest in step (II) is to describe the situation
in which we may pass from $L^p_{\varphi_1}(e\M_1e)$ to
$L^q_{\varphi_J}(e\M_1e)$ (where $1 \leq q \leq p < \infty$) by means of a
change of weights. In this regard notice that since by assumption $\varphi_2
\circ J \ll_{loc} \varphi_1$, $\varphi_2 \circ J$ is necessarily semifinite.
Given that we are only really interested in the action of $\varphi_1$ and
$\varphi_J$ on $e\M_1e$, we may assume for the sake of argument that $\varphi_2
\circ J$ is faithful.
As was noted in the preamble, care should be taken in exactly how we change
weights, if we are to end up with a composition operator. So in particular in
the noncommutative world we can not just willy nilly apply (\cite{Tp}; II.37 \&
II.38) and leave it at that. To gain some insight into what is required we take
some time to consider the semifinite case. So suppose that $\M_i$ $(i = 1, 2)$
are equipped with fns traces $\tau_1$ and $\tau_2$ respectively. From
(\cite{L1}) we see that if $J$ is in fact $\sigma$-weakly continuous (as we are
assuming here), then roughly speaking it will induce a projection preserving
bounded linear map from $L^p(\M_1, \tau_1)$ to $L^q(\M_2, \tau_2)$ if and only
if $f_J = \frac{\mathrm{d}\tau_2 \circ J}{\mathrm{d}\tau_1}$ exists as an
element of $L^r(\M_1, \tau_1)$ (where $r = \frac{p}{p-q})$ and
$$\tau_2 \circ J (a) = \tau_1(f_J^{1/2} a f_J^{1/2}) \quad
\mbox{for each} \quad a \in \M_1.$$ For any $a \in L_p(\M_1,
\tau_1) \cap M_1$ we then have
\begin{eqnarray*}
\|J(a)\|_q = (\tau_2(|J(a)|_q^q))^{1/q} &=& (\tau_2 \circ
J(|a|_q^q))^{1/q}\\ &=& (\tau_1(f_J^{1/2}|a|_q^q
f_J^{1/2}))^{1/q}\\ &\leq& \|f\|_r^{1/q} \|a\|_p
\end{eqnarray*}
(In the above $|a|_q$ denotes the so-called $q$-th symmetric
modulus discussed in \cite{L1}.) Here the first line corresponds
to the isometric embedding of $L^q(\M_1, \tau_2 \circ J)$ into
$L^q(\M_2, \tau_2)$, and the next two to the passage from
$L^p(\M_1, \tau_1)$ to $L^q(\M_1, \tau_2 \circ J)$ by means of a
change of weights. So we see that it is the derivative $f_J$ that
not only enables us to pass from $L^p(\M_1, \tau_1)$ to $L^q(\M_1,
\tau_2 \circ J)$ by means of the identity
$$\tau_2 \circ J (\cdot)
= \tau_1(f_J^{1/2} \cdot f_J^{1/2}),$$but also conditions the
boundedness of the induced map.

Passing to the general case the assumption that $J$ is normal ensures that
$\varphi_J = \varphi_2 \circ J$ is normal, in addition to being faithful and
semifinite. So for the sake of clarity we may assume for now that $\cM{1}{1} =
\cM{1}{J}$ \cite[II.37 \& II.38]{Tp}. Now let $\tr_1$ and $\tr_J$ be the
canonical trace functionals associated with $L^1_{\varphi_1}(\M_1)$ and
$L^1_{\varphi_J}(\M_1)$ respectively, and let $h_1 = \frac{\mathrm{d}
\widetilde{\varphi_1}}{\mathrm{d}\tau}$ and $h_J = \frac{\mathrm{d}
\widetilde{\varphi_J}}{\mathrm{d}\tau}$. In a simplistic world we would then by
analogy with the semifinite case hope to achieve the change of weights by means
of some positive element $f_J \in (\cM{1}{1})\widetilde{} $  for which
$\tr_J(\cdot) = \tr_1(f_J^{1/2} \cdot f_J^{1/2})$. However this is too much to
hope for in general, as the type III case is rather more exotic than the
semifinite case. This makes for a type III theory of ``composition operators''
which shows some interesting variations to the semifinite theory. If the
weights $\varphi_1$ and $\varphi_J$ actually commute, then by \cite[Corollary
VIII.3.6]{Tak} there indeed does exist some $v \geq 0$ affiliated to
$(M_1)_{\varphi_1}$ such that $$\varphi_J(\cdot) = \varphi_1(v^{1/2}\cdot v^{1/2}).$$
Although the above is already reminiscent of the
equality in the semifinite setting, it would be more useful to translate this to
a statement concerning $\tr_1$ and $\tr_J$. Now by mimicking the argument of
\cite[Proposition 2.13]{GL2} we may show that
$$\varphi_1(\sigma_{i/2}(b)c\sigma_{-i/2}(b^*)) =
\tr_1(b\mathfrak{i}^{(p)}(c)b^*) \quad b \in \mathfrak{m}_\infty, c \in
\mathfrak{n}^*.$$ Arguing formally, the fact that $v$ is affiliated to
$(M_1)_{\varphi_1}$ then seems to suggest that in the case of commuting weights
we will have $$\tr_J(h_J^{1/2}\cdot h_J^{1/2}) = \tr_1(v^{1/2} h_1^{1/2}\cdot
h_1^{1/2} v^{1/2}),$$ or in other words $$\tr_J(\mathfrak{i}_J^{(1)}(\cdot)) =
\tr_1(d\mathfrak{i}_1^{(p)}(\cdot)d^*)$$ where $d = v^{1/2}h_1^{1/(2p^*)}$. If
now $d \in L^{2p^*}_{\varphi_1}(\M_1)$, we could use H\"{o}lders's inequality
to show that then the process $\mathfrak{i}^{(p)}_1(a) =
h_1^{1/(2p)}\pi_1(a)h^{1/(2p)} \rightarrow h_J^{1/2}\pi_J(a)h_J^{1/2} =
\mathfrak{i}^{(1)}_J(a)$ $(a \in \M_1^{(0)})$ extends to a well defined bounded
map $L^p_{\varphi_1}(\M_1) \rightarrow L^1_{\varphi_J}(\M_1)$. At least for the
the case $q = 1$ the resultant map then seems to represent a means of passing
from $\mathfrak{i}^{(p)}_1(\M_1)\subset L^p_{\varphi_1}(\M_1)$ to
$\mathfrak{i}^{(1)}_J(\M_1)\subset L^1_{\varphi_J}(\M_1)$ by means of a
``change of weights'' in a way that is categorically more in line with what is
required for the construction of composition operators. Admittedly this
``change of weights'' is dependent on the manner in which $\M_1$ is embedded in
$L_p$, but this fact seems to be a challenge inherent in the type III theory.

It remains to develop a suitable strategy for dealing with the case
$L^p_{\varphi_1}(\M_1) \rightarrow L^q_{\varphi_J}(\M_1)$ where $1 < q \leq
\infty$. Formally one may consider something like
$\mathfrak{i}^{(p)}_1(a) =
h_1^{1/(2p)}\pi_1(a)h^{1/(2p)} \rightarrow h_J^{1/(2q)}\pi_J(a)h_J^{1/(2q)} =
\mathfrak{i}^{(q)}_J(a)$ $(a \in \M_1^{(0)})$. We deal with the situation by
first considering change of weights mapping acting in one specific crossed
product (say, the one given by $\varphi_1$), and then by applying the natural
isometry $\gamma$ (described in detail in \cite{Tp}; II.37 \& II.38) that
identifies this crossed product with the one given by the other weight
($\varphi_J$ in our case). The following proposition deals with the change of
weights:

\begin{proposition} \label{cw}
Let $\M$ be a von Neumann algebra with two $fns$ weights $\varphi$ and
$\varphi_0$ with $\varphi_0 \ll_{loc} \varphi$. Let $h = \frac{\mathrm{d}
\widetilde{\varphi}}{\mathrm{d}{\tau}}$ and $k = \frac{\mathrm{d}
\widetilde{\varphi_0}}{\mathrm{d}{\tau}}$. Also let $1 \leq q \leq p \leq
\infty$.

Then the following statements are equivalent:
\begin{enumerate}
\item The embedding $h^{1/(2p)}ah^{1/(2p)} \rightarrow
k^{1/(2q)}ak^{1/(2q)}$ $(a \in \M^{(0)})$ extends to a continuous map $T :
L^p_{\varphi}(\M) \rightarrow L^q_{\varphi}(\M)$;
\item for $r$ such that $\frac{1}{q} = \frac{1}{p} + \frac{1}{r}$,
there exists some $d \in L^{(2r)}_\varphi(\M) \subset \cM{}{}$, such that
$f|[dh^{1/(2p)}]|^2f = fk^{1/q}f$ for any projection $f \in \M^{(0)}$;
\item for each pair $1 \leq q_0 \leq p_0 \leq \infty$ with $(p_0, q_0)
\geq (p, q)$ (by the lexicographic ordering) and with $p_0/q_0 = p/q$, the
embedding $h^{1/(2p_0)}ah^{1/(2p_0)} \rightarrow k^{1/(2q_0)}ak^{1/(2q_0)}$ $(a
\in \M^{(0)})$ extends to a continuous linear  map $T^{(p_0,q_0)} :
L^{p_0}_{\varphi}(\M) \rightarrow L^{q_0}_{\varphi}(\M)$.
\end{enumerate}
\end{proposition}

In our construction of composition operators the operator $d$ above will then
fulfill the role played by $f_J^{1/2}$ in the semifinite setting - see the
preceding discussion.

\begin{definition}
Let $\M$ be a von Neumann algebra with an $fns$ weight $\varphi$ and let
$\varphi_0$ be a normal weight with support projection $e$ belonging to the
fixed point algebra of $\varphi$, and with $\varphi_0 \ll_{loc} \varphi$. Let
$h$ and $k$ be as in the preceding discussion. Note that our assumptions imply
that $k^{1/(2q)}ak^{1/(2q)}=k^{1/(2q)}eaek^{1/(2q)}$ is well-defined for any
$(a \in \M^{(0)})$.

Given $1 \leq p, q \leq \infty$, we say that $\M$ \emph{admits of a bounded
change of weights from $\varphi$ to $\varphi_0$ for the pair $(p, q)$}, if the
embedding $h^{1/(2p)}ah^{1/(2p)} \rightarrow k^{1/(2q)}ak^{1/(2q)}$ $(a \in
\M^{(0)})$ extends to a continuous linear map $T : L^p_{\varphi}(\M)
\rightarrow L^q_{\varphi_0}(\M )$.

Given $1 \leq r < \infty$, we say that $\M$ \emph{admits of a bounded change of
weights scale from $\varphi$ to $\varphi_0$ for the ratio $r$} if for each pair
$1 \leq q \leq p < \infty$ with $r = p/q$, the embedding $h^{1/(2p)}ah^{1/(2p)}
\rightarrow k^{1/(2q)}ak^{1/(2q)} $ extends to a continuous map $T^{(p,q)} :
L^p_\varphi(\M) \rightarrow L^q_{\varphi_0}(\M )$.

Notice that the support of $k$ is just $e$. Thus in the above definition, the
maps $T, T^{(p,q)}$ actually maps into $eL^q_{\varphi_0}(\M)e$. On canonically
identifying $eL^q_{\varphi_0}(\M)e$ with $L^q_{\varphi_0}(e\M e)$, we may therefore
equivalently speak of a bounded change of weights from $(\M, \varphi)$ to $(e\M e,
\varphi_0)$ for the pair $(p, q)$, etc.
\end{definition}

The proposition will be an easy  consequence of the following,
more general, theorem (see \cite[Theorem 2.5]{JS} and the comments
in the introduction).

\begin{theorem}\label{JSmod}
Let $\M$ be a von Neumann algebra with an $fns$ weight $\varphi$, and let $1
\leq q \leq p \leq \infty$ and  $T \in \mathrm{Hom}(L^p_\varphi(\M)_{\M},
L^q_\varphi(\M)_{\M})$ (Thus $T$ is a bounded linear map from $L^p$ to $L^q$
which is a homomorphism with respect to the right-module action of $\M$ on
$L^p$.) Then there exists $c \in L^r_\varphi(\M)$ (where $\frac{1}{q} =
\frac{1}{p} + \frac{1}{r}$) such that $T(a) = ca$ for all $a \in
L^p_\varphi(\M)$.
\end{theorem}

We first show how Proposition \ref{cw} can be deduced from the above theorem.

\begin{proof}[Proof of Proposition \ref{cw}.]

\textbf{(1) $\Rightarrow$ (2):} The implication clearly holds if $p = \infty$,
and hence we may assume that $p < \infty$.

Assume that
\begin{equation}\label{AA} h^{1/(2p)}ah^{1/(2p)}
\mapsto k^{1/(2q)}ak^{1/(2q)} \quad a \in \M^{(0)}
\end{equation} extends to a continuous map $T$ from
$L^p_{\varphi}(\M)$ to $L^q_{\varphi}(\M)$. Given any $a \in \M^{(0)}$, the
spectral resolution for selfadjoint operators ensures that we may find a
sequence of Riemann sums of the form $\sum_{i=1}^n \lambda_i e_i$ with each
$e_i$ a projection majorised by $s_r(a)$ (where $s_r(a)$ is the right support
of $a$), and $e_i$'s mutually orthogonal, which converges uniformly to $|a|^2$
in the compression $s_r(a)\M s_r(a)$. Since $\varphi(s_r(a)) < \infty$, we have
that $h^{1/(2p)}s_r(a) \in L^{2p}_{\varphi}(\M)$. Hence an application of
H\"older's inequality reveals that the terms $h^{1/(2p)}s_r(a)(\sum_{i=1}^n
\lambda_i e_i)s_r(a)h^{1/(2p)} = h^{1/(2p)}(\sum_{i=1}^n \lambda_i
e_i)h^{1/(2p)}$ must converge to $h^{1/(2p)}|a|^2h^{1/(2p)}$. Since $\varphi_0
\ll_{loc} \varphi$, we of course also have $\varphi_0(s_r(a)) < \infty$, and
hence essentially the same argument shows that the terms
$k^{1/(2q)}(\sum_{i=1}^n \lambda_i e_i)k^{1/(2q)}$ must converge to
$k^{1/(2q)}|a|^2k^{1/(2q)}$. Thus for any $a \in \M^{(0)}$ the continuity of
$T$ ensures that it will map the term $h^{1/(2p)}|a|^2h^{1/(2p)}$ onto the term
$k^{1/(2q)}|a|^2k^{1/(2q)}$. From this observation it now follows that
\begin{align*}\|[ak^{1/(2q)}]\|_{2q} &= \|k^{1/(2q)}|a|^2k^{1/(2q)}\|_q^2 =
\|T(h^{1/(2p)}|a|^2h^{1/(2p)})\|_p^2\\
&\leq \|T\|^2 \|h^{1/(2p)}|a|^2h^{1/(2p)}\|_p^2 = \|T\|^2
\|[ah^{1/(2p)}]\|_{2p}.
\end{align*}
Thus the formal map $[ah^{1/(2p)}] \mapsto [ak^{1/(2q)}]$ ($a \in \M^{(0)}$)
extends continuously to a map $T_0: L^{2p}_{\varphi}(\M) \rightarrow
L^{2q}_{\varphi}(\M)$. This map is a homomorphism with respect to the left
module action of $\M$ on $L^{2p}_{\varphi}(\M)$. Thus an application of the
{\em left} version of Theorem \ref{JSmod} now establishes (2).

\textbf{(2) $\Rightarrow$ (1):} For the converse note
that if an element $d$ of the form described in (2) exists, then given any $a
\in \M^{(0)}$, we may select a partial isometry $u$ so that
$ud^{1/2}(h^{1/(2p)}f) = k^{1/(2q)}f$ where $f = s_l(a)\vee s_r(a)$ (here
$s_l(a)$ and $s_r(a)$ are the left and right supports of $a$). A simple
application of H{\"o}lder's inequality then reveals that
\begin{eqnarray*}
\|k^{1/(2q)}ak^{1/(2q)}\|_q &=& \|udh^{1/(2p)}ah^{1/(2p)}d^*u^*\|_q\\
&\leq& \|ud\|_{(2r)} \cdot
\|h^{1/(2p)}ah^{1/(2p)}\|_p \cdot \|d^*u^*\|_{(2r)}\\
 &\leq& \||d|^2\|_{r} \cdot
\|h^{1/(2p)}ah^{1/(2p)}\|_p.
\end{eqnarray*}
Since this holds for each $a \in \M^{(0)}$, the embedding
$h^{1/(2p)}ah^{1/(2p)} \rightarrow k^{1/(2q)}ak^{1/(2q)}$ $(a \in \M^{(0)})$
therefore clearly extends to a continuous map $\widetilde{T} :
L^p_{\varphi}(\M) \rightarrow L^q_{\varphi}(\M)$.

\textbf{(1) $\Rightarrow$ (3):} Here $T^{(p,q)}$ is nothing but the unique
operator for which $T^{(p,q)} \circ \mathfrak{i}^{(p)} |{\M^{(0)}} =
\mathfrak{i}^{(q)} \circ \mathrm{id}_{\M^{(0)}}$. Now let $T^{[p,q]}$ be the
unique bounded operator on the Terp interpolation space $L^p(\M, \varphi)$ such
that $T^{[p,q]} \circ \kappa_p^{(\varphi)} = \kappa_q^{(\varphi)} \circ
T^{(p,q)}$. It is clear that $T^{[p,q]}|\kappa_p^{(\varphi)}(\M^{(0)}) =
T^{[\infty, \infty]}|\kappa_\infty^{(\varphi)}(\M^{(0)})$ where $T^{(\infty,
\infty)} = \mathrm{id}_\M$. By the reiteration property of the complex
interpolation method (\cite{BeL}; Theorem 4.6.1),
$T^{[p_0,q_0]}|\kappa_{p_0}^{(\varphi)}(\M^{(0)}) =
T^{[p,q]}|\kappa_p^{(\varphi)}(\M^{(0)})$ is bounded, which implies the
boundedness of $T^{(p_0,q_0)}$.

The implication $(3) \Rightarrow (1)$ is entirely trivial, and
hence the result follows.
\end{proof}

We begin the proof of Theorem \ref{JSmod} with two lemmas.

\begin{lemma} \label{cwl1}
Let $\M$ be a von Neumann algebra with an $fns$ weight $\varphi$ and let $1
\leq r < \infty$ be given. For any $0 < t, s < \infty$ satisfying $\frac{1}{s}
= \frac{1}{r} + \frac{1}{t}$ and any $b \in L^r_\varphi(\M)$, we have
$$\|b\|_r = \sup\{\|bg\|_s : g \in L^t_\varphi(\M), \|g\|_t \leq 1\}.$$If
$1 \leq t, s < \infty$, the formula also holds for the case $r =
\infty$.
\end{lemma}

\begin{proof}
The statement obviously holds if $b = 0$. If $b \neq 0$ we may normalise and
assume that $\|b\|_r = 1$. H\"{o}lder's inequality then ensures that $$1 \geq
\sup\{\|bg\|_s : g \in L^t_\varphi(\M), \|g\|_t \leq 1\}.$$

To see that we get equality when $1 \leq r < \infty$, consider the element of
$L^t_\varphi(\M)$ defined by $g_b = |b|^{r/t}$. Then $\|g_b\|_t =
(\|b\|_r)^{r/t} = 1$ with $\|bg_b\|_s = \tr(|bg_b|^s)^{1/s} = \tr((|b|^{1 +
(r/t)})^s)^{1/s} = \tr(|b|^r)^{1/s} = 1$ as required.

Finally let $1 \leq t = s < \infty$ and $r = \infty$. For the case $1 = t$ this
formula is known. Hence let $1 < t < \infty$. Given any $0 < \varepsilon < 1$,
we may use $L^p$ duality to select $f \in L^1_\varphi(\M)$ with $1 -
\varepsilon < \tr(bf) \leq 1$ and $\tr(|f|) = 1$. Let $u|f|$ be the polar
decomposition of $f$ and set $g_b = u|f|^{1/t} \in L^t_\varphi(\M)$. Since
$\|(|f|^{1/(t^*)})\|_{t^*} = \tr(|f|)^{1/(t^*)} = 1$, it therefore follows from
H\"{o}lder's inequality that $1 - \varepsilon < \tr(bf) =
\tr(bg_b|f|^{1/(t^*)}) \leq \|bg_b\|_t$. Now by construction $\|g_b\|_t =
\tr(|u|f|^{1/t}|^t)^{1/t} = \tr(|f|)^{1/t} =1$. Hence $1 - \varepsilon \leq
\sup\{\|bg\|_s : g \in L^t_\varphi(\M), \|g\|_t \leq 1\}$. From these
considerations it is clear that $1 = \sup\{\|bg\|_s : g \in L^t_\varphi(\M),
\|g\|_t \leq 1\}$ as required.
\end{proof}

\begin{lemma} \label{cwl2}
Let $\M$ be a von Neumann algebra equipped with an fns weight $\varphi$, and
let $h = \frac{\mathrm{d} \widetilde{\varphi}}{\mathrm{d}\tau}$. Let $2 \leq q
\leq p < \infty$ and let $T \in \mathrm{Hom}(L^p_\varphi(\M)_{\M},
L^q_\varphi(\M)_{\M})$ (Thus $T$ is a bounded linear map from $L^p$ to $L^q$
which is a homomorphism with respect to the right-module action of $\M$ on
$L^p$.) Let $e \in \M$ be a projection in $\M$ with $\varphi(e) < \infty$, and
let $d_e = T(eh^{1/p})$. Then the following holds:
\begin{enumerate}
\item For any $a \in \M$, we have that
$T([eh^{1/p}]a) = d_ea$.
\item For any $g \in L^{q^*}_\varphi(\M)$ we have $$gd_e
= T^*(g)[eh^{1/p}].$$
\item The formal map $h^{1/p}b \mapsto d_eb$ defined for all $b
\in \{a \in L^{q^*}_\varphi(\M) : \varphi(s_l(a)) < \infty\}$, extends
continuously and uniquely to a linear map $L^v_\varphi(\M) \mapsto
L^1_\varphi(\M)$ where $1 \leq v \leq \infty$ is such that $\frac{1}{v} =
\frac{1}{p} + \frac{1}{q^*}$.
\end{enumerate}
\end{lemma}

\begin{proof}
Let $T$ be a bounded linear map from $L^p_\varphi(\M)$ to $L^q_\varphi(\M)$.
By continuity and the density of $h^{1/p}\M^{(0)}$ in $L^p_\varphi(\M)$, it is
not difficult to see that (1) follows directly from the requirement that $T \in
\mathrm{Hom}(L^p_\varphi(\M)_{\M}, L^q_\varphi(\M)_{\M})$. It therefore remains
to demonstrate the validity of (2) and (3).

Next consider claim (2). For any $a \in \M$ and $g$ as in the hypothesis, we
have
\[\tr((gd_e)a) = \tr(g(d_ea)) = \tr(gT([eh^{1/p}]a)) = \tr(T^*(g)[eh^{1/p}]a).\]
It follows from this equality that $gd_e = T^*(g)[eh^{1/p}]$.

Finally consider claim (3). Given $b \in \{a \in L^{q^*}_\varphi(\M) :
\varphi(s_l(a)) < \infty\}$, it follows from claim (2) that we will then have
$$gd_eb = T^*(g)[eh^{1/p}]b, \quad g \in L^{q^*}_\varphi(\M).$$ Since $h^{1/p}b =
h^{1/p}s_l(b)b \in L^p_\varphi(\M)\cdot L^{q^*}_\varphi(\M) \subset
L^v_\varphi(\M)$, we may apply Lemma \ref{cwl1} (with $t = q^*$, $r = 1$) and
H\"{o}lder's inequality to get
\begin{eqnarray*}
\|d_eb\|_1 &=& \sup\{\|gd_eb\|_s : g \in L^{q^*}(\M), \|g\| \leq 1\}\\ &=&
\sup\{\|T^*(g)[eh^{1/p}]b\|_s : g \in L^{q^*}(\M), \|g\| \leq 1\}\\ &\leq&
\|[eh^{1/p}]b\|_v \sup\{\|T^*(g)\|_{p^*} : g \in L^{q^*}(\M), \|g\| \leq 1\}\\
&\leq& \|T\|\cdot\|[eh^{1/p}]b\|_v.
\end{eqnarray*}
Now, since $s_l(b)\in(\mathfrak{n}^{(2)})^*$, the operator $h^{1/p}s_l(b)$ is
premeasurable (in fact, even measurable). Hence, $h^{1/p}b=(h^{1/p}s_l(b))b$,
being a product of two premeasurable operators, is also premeasurable. Since
$eh^{1/p}b\subset[eh^{1/p}]b$ and $eh^{1/p}b\subset e[h^{1/p}b]$, the rigidity
of measurable operators yields $[eh^{1/p}]b=e[h^{1/p}b]$. Thus,
\[ \|T\|\cdot\|[eh^{1/p}]b\|_v\leq \|T\|.\|h^{1/p}b\|_v,\]
as required. (Here we made use of the fact that $\frac{1}{s} = 1 +
\frac{1}{q^*} = 1 - \frac{1}{p} + \frac{1}{v} = \frac{1}{p^*} + \frac{1}{v}$.)
The last part of the claim now follows from the density of $\{h^{1/p}a : a \in
L^{q^*}_\varphi(\M), \varphi(s_l(a)) < \infty\}$ in $L^v_\varphi(\M)$.
\end{proof}

\begin{proof}[Proof of Theorem \ref{JSmod}.]
As noted in \cite{JS}, the implication clearly holds if $p = \infty$,
and hence we may assume that $p < \infty$. Suppose for the sake of argument
that $p > q$. (Note: As will be seen, the proof below easily adapts for the
case $p = q$.) Notice that the above assumptions in turn ensure that $1 < r <
\infty$, and hence that $L^r$ is reflexive.

First assume that $2 \leq p < \infty$. Let $e$ be a projection in $\M$ with
$\varphi(e) < \infty$. As noted in the preceding lemma, the restriction of $T$
to $eL^p_\varphi(\M)$ is a continuous extension of the formal map $[eh^{1/p}]a
\mapsto d_ea$ ($a \in \M$) where $d_e = T([eh^{1/p}])$. Lemma \ref{cwl2} now
additionally informs us that the map $[h^{1/p}b] \mapsto d_eb$ ($b \in \{a \in
L^{q^*}_\varphi(\M) : \varphi_1(s_l(a)) < \infty\}$) is a continuous map from a
dense subspace of $L^v_\varphi(\M)$ into $L^1_\varphi(\M)$ where $v$ is such
that $\frac{1}{v} = \frac{1}{p} + \frac{1}{q^*}$. We may therefore compose this
map with the trace functional $\tr$ on $L^1_\varphi(\M)$ to get a densely
defined continuous linear functional $[h^{1/p}b] \mapsto \tr(d_eb)$ ($b \in \{a
\in L^{q^*}_\varphi(\M) : \varphi_1(s_l(a)) < \infty\}$) on $L^v_\varphi(\M)$.
Thus by $L^p$ duality there must exist $c_e \in
L^{v^*}_\varphi(\M)=L^{r}_\varphi(\M)$ with
$$\tr(c_e[h^{1/p}b]) = \tr(d_eb)$$ for all $b$ with $\varphi_1(s_l(b)) < \infty$.
It is clear that the $b$ in the formula above can be replaced with $ba$, where
$a\in \M$, and that $[h^{1/p}ba]=[h^{1/p}b]a$. This implies $c_e[h^{1/p}b]
=c_e[h^{1/p}b] = d_eb$, so that $c_e[h^{1/p}ah^{1/q^*}] = d_eah^{1/q^*}$ for
each $a\in\M^{(0)}$. Consequently, $c_eh^{1/p}ah^{1/q^*} \subset d_eah^{1/q^*}$
and the invertibility of $h$ yields $c_eh^{1/p}a \subset d_ea$. Again, by
rigidity $c_e(h^{1/p}a) = d_ea=T([eh^{1/p}]a)$ for all $a\in\M^{(0)}$.

Now let $\{e_\lambda : \lambda \in \Lambda\}$ be a mutually orthogonal family
of projections with $\varphi(e_\lambda) < \infty$ for each $\lambda$, and
$\sum_{\lambda \in \Lambda} e_\lambda = \I$. Let $a_0$ be a fixed element of
$\M^{(0)}$. For any finite subset $F$ of $\Lambda$ we have by linearity that
$$(\sum_{\lambda \in F}c_{e_{\lambda}})(h^{1/p}a_0) = T([(\sum_{\lambda
\in F}e_{\lambda})h^{1/p}]a_0).$$ The net of terms of the form $\sum_{\lambda
\in F}e_{\lambda}$ converges to $\I$ in the weak* topology, and hence the net
$\{\sum_{\lambda \in F}e_{\lambda}(h^{1/p}a_0)\}_F$ (where $F$ ranges over the
finite subsets of $\Lambda$) will converge weakly to $h^{1/p}a_0$. Thus
$T([\sum_{\lambda \in F}e_{\lambda}h^{1/p}]a_0) \rightarrow T(h^{1/p}a_0)$
weakly. If now we combine the density of $h^{1/p}\M^{(0)}$ in $L^p_\varphi(\M)$
with the previous centered equation, we get that
\begin{eqnarray*}
\|\sum_{\lambda \in F}c_{e_{\lambda}}\|_r &=& \sup\{\|(\sum_{\lambda \in
F}c_{e_{\lambda}})(h^{1/p}a)\|_q : a \in M^{(0)}, \|h^{1/p}a\|_p \leq 1\}\\
&=& \sup\{\|T((\sum_{\lambda \in F}e_{\lambda})h^{1/p}a)\|_q : a \in M^{(0)},
\|h^{1/p}a\|_p \leq 1\}\\ &\leq& \|T\|,
\end{eqnarray*}
since again $[eh^{1/p}]a=e[h^{1/p}a]$ for $e,a\in\M^{(0)}$. Therefore by the
weak compactness of the unit ball of $L^r$ we may select a subnet of terms of
the form $\sum_{\lambda \in \widetilde{F}}c_{e_{\lambda}} \in L^r$ (where
$\widetilde{F} \subset \Lambda$ is finite) converging to some $c \in L^r$. (In
the case $p = q$ we would have $r = \infty$. Hence we could then use weak*
compactness instead of weak compactness.) By now taking limits it follows that
$$c(h^{1/p}a_0) = \lim_{\widetilde{F}}(\sum_{\lambda \in
\widetilde{F}}c_{e_{\lambda}})(h^{1/p}a_0) =
\lim_{\widetilde{F}}T((\sum_{\lambda \in
\widetilde{F}}e_{\lambda})(h^{1/p}a_0)) = T(h^{1/p}a_0).$$Since $a_0$ was an
arbitrary element of $\M^{(0)}$, we may now finally appeal to the density of
$h^{1/p}\M^{(0)}$ in $L^p$, to conclude that as required $$cb = T(b) \quad
\mbox{for all} \quad b \in L^p.$$

Notice that everything we have done so far is entirely
symmetrical, and hence we may similarly prove that if $2 \leq p
\leq \infty$, then all left $\M$-module homomorphisms from
$L^p(\M)$ to $L^q(\M)$ are right multiplication operators induced
by some $c \in L^r(\M)$.

Now suppose that $1 \leq p < 2$. It is an exercise to show that $T
: L^p(\M) \to L^q(\M)$ is a right $\M$-module homomorphism if and
only if $T^* : L^{q^*}(\M) \to L^{p^*}(\M)$ is a left $\M$-module
homomorphism, and that $T$ is a left multiplication operator
induced by some element $c \in L^r(\M)$ if and only if $T^*$ is a
right multiplication operator induced by the same element $c$
(notice that here $\frac{1}{p^*} = \frac{1}{q^*} + \frac{1}{r})$.
Notice for example that if for some $a \in \M$ we have that $T(b)a
= T(ba)$ for every $b \in L^p(M)$, we will then have that
$tr(T^*(ax)b) = tr(axT(b)) = tr(xT(b)a) = tr(xT(ba)) =
tr(T^*(x)ba) = tr(aT^*(x)b)$ for every $b \in L^p(\M)$ and every
$x \in L^{q^*}(\M)$. Thus we then clearly have that $T^*(ax) =
aT^*(x)$ for every $x \in L^{q^*}(\M)$. Therefore since $1 \leq p
< 2$ forces $1 \leq q < 2$ (or equivalently $2 < q^* \leq
\infty$), the present case clearly follows by duality from the
case $2 \leq p \leq \infty$.
\end{proof}

\subsubsection*{Step (III) : Applying the Jordan morphism}

We start with the simplest case when $\mathcal{B}=\M_2$ and $J:\M_1 \to \M_2$
is a Jordan $*$-isomorphism of $\M_1$ onto $\M_2$. The challenge is then to
find a natural canonical way of isometrically identifying
$L^q_{\varphi_J}(\M_1)$ (where $\varphi_J = \varphi_2 \circ J$) with
$L^q_{\varphi_2}(\M_2)$. In a sequence of papers (\cite{W1} - \cite{W3})
Keiichi Watanabe developed just such a construction. (See for example \S 3 of
\cite{W3} and the discussion preceding 3.1 of \cite{W5}.) All we need to do is
to apply Watanabe's construction to $J^{-1}$ to get the following

\begin{lemma}\label{wat} Let $J$ be a bijective Jordan $*$-isomorphism.
Then $J$ canonically extends to a Jordan $*$-isomorphism $\widetilde{J}$ from
$(\M_1 \rtimes_{\sigma^{J}}{\mathbb R})\widetilde{} $onto $(\M_2
\rtimes_{\sigma^{2}}{\mathbb R})\widetilde{}$ which canonically identifies $h_J
= \frac{\mathrm{d}\widetilde{\varphi_J}}{\mathrm{d}\tau_J}$ with $h_2 =
\frac{\mathrm{d}\widetilde{\varphi_2}}{\mathrm{d}\tau_2}$, and isometrically
identifies $L^p_{\varphi_J}(\M_1)$ with $L^p_{\varphi_2}(\M_2)$. (Here
$\varphi_J = \varphi_2 \circ J$.)
\end{lemma}

(Note that in the computation in the middle of p 275 of \cite{W1} it is shown
that $\widetilde{J}$ takes the \emph{shift} map $\lambda^1_s$ onto
$\lambda^2_s$. This fact together with the continuity of $\widetilde{J}$ in the
topology of convergence in measure, now ensures that speaking loosely $h_J$
\emph{maps onto} $h_2$ with respect to this identification.)

We now move on to the more general case when the image of $\M_1$ under $J$ is
not necessarily a von Neumann algebra. Let $e$ be the support projection of
$\varphi_J$. We remind the reader that $e$ belongs to the center of $\M_1$. Let
$z$ be a central projection in $\mathcal{B}$ such that $zJ$ is a
*-homomorphism and $(1-z)J$ is a
*-antihomomorphism. Since the kernels of $zJ$ and $(\I-z)J$ are both two-sided
ideals in $\M_1 e$, there exist central projections $e_z$ and $e_{\I-z}$ in
$\M_1 e$ such that $\ker(zJ)=\M_1(\I-e_z)$ and $\ker((\I-z)J)=\M_1(\I-e_{\I-z})$.
Note that $e_z$ is the support of $\varphi_z = \varphi_2\circ zJ$ and $e_{\I-z}$ is the
support of $\varphi_{\I-z} = \varphi_2\circ (\I-z)J$. Note also that $zJ(\M_1e_z) =
\mathcal{B}z$ and $(\I-z)J(\M_1e_{\I-z})=\mathcal{B}(\I-z)$. This follows easily from
the fact that the smallest von Neumann algebra containing $J(\M_1)$ must also contain the
projection $z$, and by then realizing that the direct sum of $zJ(\M_1e_z)$ and
$(\I-z)J(\M_1e_{\I-z})$ is a von Neumann algebra contained in $\mathcal{B}$, and
containing both $z$ and $J(\M_1)=J(\M_1 e)$ (obviously $e=e_z\vee e_{\I-z}$).
Therefore $\mathcal{B} = zJ(\M_1 e_z) \oplus (\I-z)J(\M_1e_{1-z})$  Thus $zJ$
restricted to $\M_1e_z$ is a *-isomorphism of $\M_1e_z$ onto
$\mathcal{B}z$ and $(\I-z)J$ restricted to $\M_1e_{\I-z}$ is a *-antiisomorphism
of $\M_1e_{\I-z}$ onto $\mathcal{B}(\I-z)$, and Lemma \ref{wat} shows that the
spaces $L^q_{(\varphi_z , \varphi_{\I-z})}(\M_1 e_z \times \M_1 e_{\I-z})$
and $L^q_{z\varphi_\mathcal{B} \oplus (\I-z)\varphi_\mathcal{B} }(\mathcal{B} z
\oplus \mathcal{B} (\I-z))=L^q_{\varphi_\mathcal{B} }(\mathcal{B})$ are
isometric. The `direct product' notation for the first space is used remind the
reader that $e_z$ and $e_{1-z}$ are not, in general, orthogonal to each other.
We denote the isometry mentioned above by $W_J$. With reference to Lemma \ref{wat},
it is clear that this isometry is constructed from the action of $(zJ, (\I-z)J)$
on $(\M_1 e_z \times \M_1 e_{\I-z})$. Setting $h_z =
\frac{\mathrm{d}\widetilde{\varphi_z}}{\mathrm{d}\tau_z}$ and
$h_{(\I-z)} = \frac{\mathrm{d}\widetilde{\varphi_{(\I-z)}}}{\mathrm{d}\tau_{\I-z}}$,
it is therefore an exercise to see that  $W_J$ will map elements of the form
$(h_z^{1/(2q)}\pi_{z}(a)h_z^{1/(2q)}, h_{\I-z}^{1/(2q)}\pi_{\I-z}(b)h_{\I-z}^{1/(2q)})$
(where $a, b \in \M^{(0)}$), onto $h_2^{1/(2q)}\pi_{\mathcal B}(zJ(a)+(\I-z)J(b))h_2^{1/(2q)}$.

\subsubsection*{Step (IV) and (V): Passing from
$L^q_{\varphi_2}(\mathcal{B})$ to $L^q_{\varphi_2}(\M_2)$}

Let us remind the reader that $\mathcal{B}$ is the von Neumann subalgebra of
$\M_2$ generated by $J(\M_1)$. We shall need the following results:

\begin{proposition}\label{JXprop}
Let $(\M_i,\varphi_i), i=1,2$ be von Neumann algebras with \emph{fns} weights,
and let $j:\M_1\to\M_2$ be a positive map satisfying $\varphi_2\circ j\leq C
\varphi_1$ for some $C>0$. Denote by $j^{(p)}, 1\leq p<\infty$ the maps defined
on a dense subspace of $L^p_{\varphi_1}(\M_1)$ by $h_1^{(1/2p)}a
h^{(1/2p)}\mapsto h_2^{(1/2p)}j(a) h_2^{(1/2p)}$, where $h_i$ are the densities
of the dual weights of $\varphi_i$  with respect to canonical traces on the
corresponding crossed products. Then all the $j^{(p)}$'s extend to bounded
linear operators from $L^p_{\varphi_1}(\M_1)$ into $L^p_{\varphi_2}(\M_2)$.
\end{proposition}

The proposition was proved for finite weights by Junge and Xu
\cite[Theorem 5.1]{JX}, where the norms of the mappings are also
calculated. The proof is based on Haagerup's lemma \cite[Lemma
1.1]{H} (see also \cite[Lemma VII.1.9]{Tak}). The lemma shows
essentially that for a self-adjoint element $a$ of
$\mathfrak{m}_\varphi$,
$$\|h_1^{1/2}ah_1^{1/2}\|=\inf\{\varphi(b)+\varphi(c):a=b-c,
b,c\in\mathfrak{m}_\varphi^+\}.$$ Note that the assumed inequality
gives boundedness of our mappings on positive elements, and
Haagerup's lemma allows us to extend the bound to self-adjoint
elements. That this implies boundedness of the mappings on arbitrary
elements is trivial. Since the lemma is true for weights, the
proposition is also true for weights, essentially without changes.

\begin{lemma} Let $\M_i$ and $\varphi_i$ $(i = 1, 2)$ be as
before and let $J: \M_1 \rightarrow \M_2$ be a normal Jordan
$*$-isomorphism. If $\varphi_2 \circ J \ll_{loc} \varphi_1$, then
$\varphi_2$ is semifinite on $\mathcal{B}$.
\end{lemma}

\begin{proof} It is enough to show that for any projection $e\in
\mathcal{B}$ there exists a projection $f\in\mathcal{B}$ such that
$f\leq e$ and $\varphi_2(f)<\infty$. Let $z$ be the central
projection in $\mathcal{B}$ such that $a\mapsto zJ(a)$ is a
*-homomorphism and $a\mapsto (J(\I)-z)J(a)$ is a
*-antihomomorphism. Note that both $zJ(\M_1)$ and
$(J(\I)-z)J(\M_1)$ are von Neumann algebras, as the
(anti)homomorphic images of a von Neumann algebra, with both the
homomorphism and the antihomomorphism normal (see
\cite[Proposition III.3.12]{Tak}). As noted at the close of the
discussion pertaining to Step (III), the direct sum of these two
von Neumann algebras is precisely $\mathcal{B}$. Thus there exist
projections $e_1$ and $e_2$ in
$\M_1$ such that $e=zJ(e_1)+(J(\I)-z)J(e_2)$. Choose now
projections $f_1, f_2$ in $\M_1$ so that $f_i\leq e_i$ and
$\varphi_1(f_i)<\infty$ for $i=1,2$. Put
$f=zJ(f_1)+(J(\I)-z)J(f_2)$. Then $f$ is a projection in
$\mathcal{B}$, $f\leq e$ and $\varphi_2(f)<\infty$.
\end{proof}

Let us apply the proposition to the natural embedding $j$ of the von
Neumann algebra $\mathcal{B}$ with weight $\varphi_2$ restricted to
the algebra, into the algebra $\M_2$. The inequality required for
the lemma is clearly satisfied with constant $1$.

\begin{remark} \label{IVandV}
The maps $j^{(p)}$ are especially simple if the algebra $\mathcal{B}$ is
invariant under the modular group for the couple $(\M_2, \varphi_2)$. Then the
space $L^p_{\varphi_2|\mathcal{B}}(\mathcal{B})$ can be treated as a subspace
of $L^p_{\varphi_2}(\M_2)$ and $j^{(p)}$ is the natural embedding. To see this,
we can mimic the proof of Lemma \ref{fixed}. In fact up to this canonical
embedding, the maps $\mathfrak{i}^{(p)} \circ J$ will in this case yield
essentially identical terms on $\M_1^{(0)}$ for both $\M_2$ and $\mathcal{B}$.
Thus in dealing with composition operators we may then freely replace $\M_2$
with $\mathcal{B}$. To see this note that in this case $J(\I)$ (the unit of
$\mathcal{B}$) will be a fixed point of the modular group generated by
$\varphi_2$. In this regard observe that the identity
$\sigma_t^{\varphi_2}(J(\I)a) = \sigma_t^{\varphi_2}(a)$ for all $a \in
\mathcal{B}$ and all $t \in \mathbb{R}$, ensures that
$\sigma_t^{\varphi_2}(J(\I))$ is an identity for
$\sigma_t^{\varphi_2}(\mathcal{B}) = \mathcal{B}$, and hence that
$\sigma_t^{\varphi_2}(J(\I)) = J(\I)$ for all $t \in \mathbb{R}$. By Lemma
\ref{fixed}, the density of $\varphi_2$ restricted to $J(\I)\M_2J(\I)$ may then
be identified with $J(\I)h_2 = h_2J(\I)$ (where as before $h =
\frac{\mathrm{d}\widetilde{\varphi_2}}{\mathrm{d}\tau_2}$). In addition by
\cite[IX.4.2]{Tak} there exists a faithful normal conditional expectation $E:
J(\I)\M_2J(\I) \rightarrow \mathcal{B}$ such that $\varphi_2 \circ E =
\varphi_2$ on $J(\I)\M_2J(\I)$. Hence \cite[4.8]{G} assures us that the density
of the restriction of $\varphi_2$ to $\mathcal{B}$, may be identified with that
of the restriction to $J(\I)\M_2J(\I)$, described above. Thus up to canonical
inclusion we have $\mathfrak{i}^{(p)} \circ J(a) =
h_2^{1/(2p)}J(a)h_2^{1/(2p)}$ $(a \in \M_1^{(0)})$ for both $\M_2$ and
$\mathcal{B}$.
\end{remark}

\subsubsection*{The main result.} With a description of steps (I)
-- (V) now finally behind us, we are ready to give a description of a large
class of Jordan $*$-morphisms which do yield composition operators.

\begin{lemma}
Let $\M_i$ ($i = 1, 2$) be von Neumann algebras equipped with faithful normal
weights $\varphi_i$, and let $J: \M_1 \to \M_2$ be a normal $*$-(anti)homomorphism.
(In this case $\mathcal{B} = J(\M_1)$.)
Denote the support projection of $\varphi_2 \circ J = \varphi_J$ by $e$.
Suppose that $\sigma^{\varphi_2}_t(\mathcal{B}) = \mathcal{B}$ for each $t \in
\mathbb{R}$. Then for each $1 \leq q \leq p < \infty$,  $J$ canonically induces
a composition operator from $L^p_{\varphi_1}(\M_1)$ to $L^q_{\varphi_2}(\M_2)$
if and only if firstly $\varphi_J \ll_{loc} \varphi_1$, and secondly $\M_1$
admits of a bounded change of weights from $\varphi_1$ to $\varphi_J$ for the
pair $(p, q)$.
\end{lemma}

\begin{proof} By steps (IV) and (V), the assumption that
$\sigma^{\varphi_2}_t(\mathcal{B}) = \mathcal{B}$ for each $t \in \mathbb{R}$,
enables us to reduce to the case where $J$ is surjective (see Remark
\ref{IVandV}). The rest of the proof is then essentially contained in step (I),
step (III), and Proposition \ref{cw}.
\end{proof}

Before actually extracting our main theorem from the above lemma, we need one
final technical observation regarding commuting weights. The result is surely reflected
in the literature somewhere, but we have been unable to find a reference, and hence
elect to prove the relevant lemmas in full.

Given two densely defined closed operators affiliated to some von Neumann algebra $\M$, we say that such operators \emph{commute} if they are affiliated to a common abelian von Neumann subalgebra of $\M$ (or equivalently if their spectral projections commute).

\begin{lemma} Let $\M$ be a von Neumann algebra equipped with faithful normal semifinite weight $\varphi$, and let $d\in \M$. Then $d$ commutes with $h_\varphi = \frac{\mathrm{d}\widetilde{\varphi}}{\mathrm{d}\tau}$ if and only if $d\in\M^\varphi$, the centralizer of $\varphi$ in $\M$.
\end{lemma}
\begin{proof}
For a faithful normal semifinite  weight $\psi$ on $\M$ and a positive self-adjoint densely defined operator $h$ affiliated with $\M^\psi$, the weight $\psi_h$ is defined as in \cite[Lemma VIII.2.8]{Tak}.

We can assume that the crossed product is built using the weight $\varphi$. By definition of the Radon-Nikodym derivative, $\tilde{\varphi}=\tau_{h_\varphi}$.
Thus by formula (11) in chapter II of \cite{Tp} and \cite[Lemma VIII.2.10]{Tak},
\[
\sigma_t^\varphi(d)=\sigma_t^{\tilde{\varphi}}(d)=h_\varphi^{it}dh_\varphi^{-it}.
\]
From this it clearly follows that $d\in\M^\varphi$ if and only if $h_\varphi^{it}dh_\varphi^{-it} = d$. Since the latter equality holds if and only if $d$ and $h_\varphi$ commute (see \cite[Theorem 1.VIII.13]{RS}), we are done.
\end{proof}

\begin{lemma}
Let $\varphi,\psi$ be faithful normal semifinite weights on $\M$. Then $\varphi$ and $\psi$ commute (in the sense of satisfying the conditions in \cite{Tak}, Corollary VIII.3.6) if and only if the densities $h_\varphi$ and $h_\psi$ commute.
\end{lemma}

\begin{proof}
If $\varphi$ and $\psi$ commute, then there exists a nonsingular positive self-adjoint densely defined operator $d$ affiliated with the algebra $\M^\varphi$ such that $\psi=\varphi_d$. Using formula (12) from chapter II of \cite{Tp}, \cite[4.8]{Str} and the chain rule for the Connes cocycle derivative, we conclude that
\[
h_\psi^{it}=(D\tilde{\varphi_d}:D\tau)_t=(D\tilde{\varphi_d}:D\tilde{\varphi})_t(D\tilde{\varphi}:D\tau)_t=
(D\varphi_d:D\varphi)_th_\varphi^{it}=d^{it}h_\varphi^{it}.
\]
Since $d$ is affiliated with the algebra $\M^\varphi$, $d^{it}$ must commute with $h_\varphi$ by the previous lemma. Hence,
\[h_\psi^{it}h_\varphi^{is}=d^{it}h_\varphi^{it}h_\varphi^{is}=h_\varphi^{is}d^{it}h_\varphi^{it}=
h_\varphi^{is}h_\psi^{it},
\]
which means, again by \cite[Theorem 1.VIII.13]{RS}, that $h_\varphi$ and $h_\psi$ commute.

Conversely, assume that $h_\varphi$ and $h_\psi$ commute. Let $\mathcal{A}$ be the abelian von Neumann algebra generated by the two operators. Since $h_\varphi^{-1}$ is a densely defined positive self-adjoint operator affiliated with $\mathcal{A}$, we can put $d=h_\psi\cdot h_\varphi^{-1}$ (see \cite{KR}, Theorem 5.6.15 (iii)). Obviously, $d$ commutes with both $h_\varphi$ and $h_\psi$, and, for each $t\in\mathbb{R}$, $d^{it}=h_\psi^{it}h_\varphi^{-it}$ (for if $\mathcal{A}$ is identified with the algebra of continuous functions on an extremely disconnected compact Hausdorff space $X$, then $f_d^{it}=f_\psi^{it}f_\varphi^{-it}$ for functions $f_d, f_\psi$ and $f_\varphi$ corresponding to the operators $d, h_\psi$ and $h_\varphi$, respectively). Moreover,
\[\sigma_s^\varphi(d^{it})=h_\varphi^{is}d^{it}h_\varphi^{-is}=d^{it},\]
so that $d$ is affiliated with $\M^\varphi$. Hence (using formula (11) from chapter II of \cite{Tp}, \cite[4.8]{Str} and the chain rule for the Connes cocycle derivative),
\[
(D\psi:D\varphi)_t=(D\tilde{\psi}:D\tilde{\varphi})_t=
(D\tilde{\psi}:D\tau)_t(D\tilde{\varphi}:D\tau)_t^*=h_\psi^{it}h_\varphi^{-it}\in\M^\varphi,
\]
which guarantees, by the Pedersen-Takesaki theorem (see \cite[4.10(iii)]{Str}) that $\varphi$ and $\psi$ commute.
\end{proof}

\begin{definition}
We say that two normal semifinite weights on $\M$ commute if the support projections of the weights commute, and the restrictions of the weights to the product of support projections also commute.
\end{definition}

\begin{theorem}\label{main}
Let $\M_i$ ($i = 1, 2$) be as before, and let $1 \leq r < \infty$ and $1 \leq q
\leq p \leq \infty$ be given. Further, let $J:\M_1 \to \M_2$ be a normal Jordan
$*$-morphism. Finally let $z$ be a central projection in $\mathcal{B}$ for which
$zJ$ and $(\I-z)J$ are respectively a $*$-homomorphism and a $*$-antihomomorphism.
As in Step (III) we write $e_z, e_{\I-z}$ for the central support projections of
$zJ$ and $(\I-z)J$, and set $\varphi_z = \varphi_2 \circ zJ$, $\varphi_{\I-z} =
\varphi_2 \circ (\I-z)J$.
\begin{enumerate}

\item[(a)] For each
pair $(p, q)$ with $p/q = r$, $J$ canonically induces a composition operator
from $L^p_{\varphi_1}(\M_1)$ to $L^q_{\varphi_2}(\M_2)$ if and only if
$\varphi_J \ll_{loc} \varphi_1$ and $\M_1$ admits of a bounded change of
weights scale from $\varphi_1$ to $\varphi_J$ for the ratio $r$.

\item[(b)] Consider the following statements:
\begin{enumerate}
\item[(1)] For $1 \leq q \leq p \leq \infty$, $J$ canonically induces a
composition operator from $L^p_{\varphi_1}(\M_1)$ to $L^q_{\varphi_2}(\M_2)$.
\item[(2)] $\varphi_J \ll_{loc} \varphi_1$ and $\M \times \M$ admits of a bounded
change of weights from $(\varphi_1, \varphi)$ to $(\varphi_z, \varphi_{\I-z})$
for the pair $(p, q)$.
\item[(3)] $\varphi_J \ll_{loc} \varphi_1$ and $\M$ admits
of a bounded change of weights from $\varphi_1$ to $\varphi_J$ for the pair
$(p, q)$.
\end{enumerate}
In general $b(3) \Rightarrow b(2) \Rightarrow b(1)$. If
$\sigma^{\varphi_2}_t(\mathcal{B}) = \mathcal{B}$ for each $t \in \mathbb{R}$,
then statements $b(1)$ and $b(2)$ are equivalent. If the weights $\varphi_z$ and $\varphi_{\I-z}$
commute, $b(2)$ and $b(3)$ are equivalent. If $\mathcal{B} = J(\M_1)$, $z$ can be
chosen so that $\varphi_z$ and $\varphi_{\I-z}$ commute.
\end{enumerate}
\end{theorem}

\begin{proof}
Throughout the proof we will let $\mathcal{B}$ and $e$ be as before. As noted
in steps (I) and (II), the centrality of $e$ enables us to assume that
$\varphi_2 \circ J$ is faithful (ie. that $e = \I$). Again for the sake of
simplicity we will now suppress the technicalities inherent in \cite[II.37 \&
II.38]{Tp}, and identify the crossed products of $\M_1$ with $\varphi_1$, and
$\M_1$ with $\varphi_J$.

\medskip
\noindent\textbf{(a):} To see the \emph{only if} part, assume that for each
pair $(p, q)$ with $p/q = r$, $J$ canonically induces a composition operator
from $L^p_{\varphi_1}(\M_1)$ to $L^q_{\varphi_2}(\M_2)$. Then for $q = 1$, $p =
r$, the map $h_1^{1/(2r)}ah_1^{1/(2r)} \mapsto h_2^{1/2}J(a)h_2^{1/2}$ $(a \in
\M_1^{(0)}$) extends to a bounded map $C_J : L^r_{\varphi_1}(\M_1) \to
L^1_{\varphi_2}(\M_2)$. (Here we have as before that $h_i = \frac{\mathrm{d}
\widetilde{\varphi_i}}{\mathrm{d}\tau_i}$.)

On composing the operator $C_J$ with the trace functional $tr_2$ on
$L^1_{\varphi_2}(\M_2)$, we obtain a positive bounded linear functional $tr_2
\circ C_J : L^r_{\varphi_1}(\M_1) \to \mathbb{C}$. Hence by $L^p$ duality there
exists $b \in L^{r^*}_{\varphi_1}(\M_1)^+$ with $tr_1(bc) = tr_2(C_J(c))$ for
each $c \in L^r_{\varphi_1}(\M_1)$. Thus
$$tr_1(b(h_1^{1/(2r)}ah_1^{1/(2r)})) = tr_2(C_J(h_1^{1/(2r)}ah_1^{1/(2r)}))
= tr_2(h_2^{1/2}J(a)h_2^{1/2})$$ for each $a \in \M_1^{(0)}$. But with $e$ as
in the hypothesis and $k = \frac{\mathrm{d} \widetilde{\varphi_2\circ
J}}{\mathrm{d}\tau_1}$, we have by \cite[2.13(a)]{GL2} that
$$tr_2(h_2^{1/2}J(a)h_2^{1/2}) = \varphi_2(J(a)) = tr_1(k^{1/2}ak^{1/2}) \quad
\mbox{for all } a \in \M_1^{(0)}.$$ But then
$$tr_1(b(h_1^{1/(2r)}ah_1^{1/(2r)})) =
tr_1(k^{1/2}ak^{1/2}) \quad \mbox{for all } a \in \M^{(0)}_1.$$ This suffices
to force $[fh_1^{1/(2r)}]b(h^{1/(2r)}f) = fkf$ for any projection $f$ with
$\varphi_1(f) < \infty$. The claim follows.

For the \emph{if} part suppose that $\varphi_J \ll_{loc} \varphi_1$ and $\M_1$
admits of a bounded change of weights scale from $(\M_1, \varphi_1)$ to
$(e\M_1e, \varphi_J)$ for the ratio $r$. A perusal of steps (I) to (III) will
reveal that this is enough to ensure that for each $1 \leq p, q \leq \infty$
with $\frac{p}{q} = r$, $J$ induces a (generalised) composition operator from
$L^p_{\varphi_1}(\M_1)$ to $L^q_{\varphi_2|_\mathcal{B}}(\mathcal{B})$ (for
further details see the proof of (b) below). An application of Proposition
\ref{JXprop} to the injection $\mathcal{B} \rightarrow \M_2$ now completes the
proof.

\medskip

\noindent{\textbf(b)} To facilitate the task of reading the proof we write
below the explicit decomposition of the composition operator $C_J$ for the
normal Jordan *-morphism $J$.
\begin{multline*}
h_1^{1/(2p)}\pi_1(a) h_1^{1/(2p)} \mapsto h_1^{1/(2p)}\pi_1(eae) h_1^{1/(2p)}\\
\mapsto k^{1/(2q)}\pi_1(eae) k^{1/(2q)}  \mapsto h_J^{1/(2q)}\pi_J(eae) h_J^{1/(2q)}\\
 \mapsto
(h_z^{1/(2q)}\pi_{zJ}(eae)h_z^{1/(2q)}, h_{\I-z}^{1/(2q)}\pi_{(\I-z)J}(eae)h_{\I-z}^{1/(2q)})\\
\mapsto h_{\mathcal{B}}^{1/(2q)}\pi_{\mathcal{B}}(J(a)) h_{\mathcal{B}}^{1/(2q)}
\mapsto h_2^{1/(2q)}\pi_2 (J(a)) h_2^{1/(2q)}
\end{multline*}

In the above scheme $h_z = \frac{\mathrm{d}\widetilde{\varphi_z}}{\mathrm{d}\tau_z}$ and
$h_{(\I-z)} = \frac{\mathrm{d}\widetilde{\varphi_{(\I-z)}}}{\mathrm{d}\tau_{\I-z}}$. In
addition to the simplifying assumptions made at the start of the proof, we will
in $b$ also use \cite[II.37]{Tp} to identify the crossed product of $\M_1e_z$ and
$\varphi_z$, with $e_z(\cM{1}{J})$. Similarly the crossed product of $\M_1e_{\I-z}$ and
$\varphi_{\I-z}$, is identified with $e_{\I-z}(\cM{1}{J})$. All of these simplifying
assumptions have the effect of identifying $k$ with $h_J$, and of forcing $h_J =
h_z + h_{\I-z}$. (This last equality is a simple consequence of the fact that
$\varphi_J = \varphi_z + \varphi_{\I-z}$.)

\medskip

\noindent{\textbf{b(3) $\Rightarrow$ b(2)}}:
Suppose that $\varphi_J \ll_{loc} \varphi_1$. Since $\varphi_J \geq \varphi_z,
\varphi_{\I-z}$, the continuity of the maps
$$h_J^{1/(2q)}ah_J^{1/(2q)} \mapsto h_z^{1/(2q)}ae_zh_z^{1/(2q)}$$ and
$$h_J^{1/(2q)}ah_J^{1/(2q)} \mapsto h_{\I-z}^{1/(2q)}ae_{\I-z}h_{\I-z}^{1/(2q)})$$(where
$a \in \M_1^{(0)}$), is an easy consequence of Proposition \ref{JXprop}. If therefore
$b(3)$ holds, we merely need to compose the above maps with the given bounded change of weights
from $\varphi_1$ to $\varphi_J$ for the pair $(p, q)$, to see that
$$(h_1^{1/(2p)}ah_1^{1/(2p)}, h_1^{1/(2p)}bh_1^{1/(2p)}) \mapsto
(h_z^{1/(2q)}ae_zh_z^{1/(2q)}, h_{\I-z}^{1/(2q)}be_{\I-z}h_{\I-z}^{1/(2q)})$$(where
$a, b \in \M_1^{(0)}$) is continuous.

\medskip

\noindent{\textbf{b(2) $\Rightarrow$ b(1)}}: Suppose that $b(2)$ holds. We show
that the hypothesis of $b(2)$ is strong enough to ensure that $J$ induces a
composition operator from $L^p_{\varphi_1}(\M_1)$ to
$L^q_{\varphi_2}(\mathcal{B})$. The conclusion will then follow from applying
Proposition \ref{JXprop} to the inclusion $\mathcal{B} \rightarrow \M_2$. In
the remainder of the proof of this implication, we may therefore assume that
$\mathcal{B} = \M_2$. From Lemma \ref{fixed} and the discussion following Lemma
\ref{supp}, it is clear that $L^q_{\varphi_2}(\M_2) = zL^q_{\varphi_2}(\M_2)
\oplus (\I - z)L^q_{\varphi_2}(\M_2) = L^q_{\varphi_2}(z\M_2) \oplus
L^q_{\varphi_2}((\I-z)\M_2)$. Taking into account the action of the isometry
described in the discussion following Lemma \ref{wat}, it is clear that $J$ will
induce the required composition operator from $L^p_{\varphi_1}(\M_1)$ to
$L^q_{\varphi_2}(\M_2)$, if the map $$h_1^{1/(2p)}a h_1^{1/(2p)} \mapsto
(h_z^{1/(2q)}ae_zh_z^{1/(2q)}, h_{\I-z}^{1/(2q)}ae_{\I-z}h_{\I-z}^{1/(2q)})$$(where
$a \in \M_1^{(0)}$) extends to a bounded linear map from $L^p_{\varphi_1}(\M_1)$ to
\newline $L^q_{(\varphi_z , \varphi_{\I-z})}(\M_1 e_z \times \M_1 e_{\I-z})$. But this map is
just the given bounded change of weights described in $b(2)$, composed with the bounded
injection $$L^p_{\varphi_1}(\M_1) \to L^p_{(\varphi_1, \varphi_1)}(\M_1 \times \M_1):
a \mapsto (a, a).$$Hence the claim follows.

\medskip

\noindent{\textbf{b(1) $\Rightarrow$ b(2)}}: Suppose that $\varphi_J \ll_{loc}
\varphi_1$. Assume that $\sigma^{\varphi_2}_t(\mathcal{B}) = \mathcal{B}$ for
each $t \in \mathbb{R}$. Steps (IV) and (V) ensure that we may then assume
$\M_2 = \mathcal{B}$ (see Remark \ref{IVandV}). Under this assumption we may
therefore select a central projection $z \in \M_2$ so that $zJ$ is a
$*$-homomorphism onto $z\M_2$ and$(\I - z)J$ a $*$-antihomomorphism onto $(\I - z)\M_2$.

Now suppose that $b(1)$ holds. From the action of the isometry described in
the discussion following Lemma \ref{wat}, it is clear that $b(1)$ is exactly
equivalent to the continuity of the map $$h_1^{1/(2p)}ah_1^{1/(2p)} \mapsto
(h_z^{1/(2q)}ae_zh_z^{1/(2q)}, h_{\I-z}^{1/(2q)}ae_{\I-z}h_{\I-z}^{1/(2q)})$$(where
$a \in \M_1^{(0)}$). Thus each of $h_1^{1/(2p)}ah_1^{1/(2p)} \mapsto
h_z^{1/(2q)}ae_zh_z^{1/(2q)}$ and \newline $h_1^{1/(2p)}ah_1^{1/(2p)} \mapsto
h_{\I-z}^{1/(2q)}ae_{\I-z}h_{\I-z}^{1/(2q)}$ are separately continuous, which in turn
is sufficient to force the continuity of $$(h_1^{1/(2p)}ah_1^{1/(2p)},
h_1^{1/(2p)}bh_1^{1/(2p)}) \mapsto (h_z^{1/(2q)}ae_zh_z^{1/(2q)},
h_{\I-z}^{1/(2q)}be_{\I-z}h_{\I-z}^{1/(2q)})$$(where $a, b \in \M_1^{(0)}$) as
required.

\medskip

\noindent{\textbf{b(2) $\Rightarrow$ b(3)}}: Suppose that $b(2)$ holds, and assume
that $\varphi_z, \varphi_{\I-z}$ commute. By the lemma, this has the effect of ensuring that $h_z$
and $h_{\I-z}$ are commuting affiliated operators. On composing the map
$$h_1^{1/(2p)}ah_1^{1/(2p)} \mapsto (h_1^{1/(2p)}ah_1^{1/(2p)}, h_1^{1/(2p)}ah_1^{1/(2p)})
\quad a \in \M_1^{(0)}$$with the given change of weights, we obtain the continuity of the map $$h_1^{1/(2p)}ah_1^{1/(2p)} \mapsto (h_z^{1/(2q)}ae_zh_z^{1/(2q)},
h_{\I-z}^{1/(2q)}ae_{\I-z} h_{\I-z}^{1/(2q)}) \quad a \in \M_1^{(0)}.$$
By continuity this map will also map terms of the form $h_1^{1/(2p)}|a|^2h_1^{1/(2p)}$
(where $a \in \M_1^{(0)}$), onto the terms $(h_z^{1/(2q)}|a|^2e_zh_z^{1/(2q)},
h_{\I-z}^{1/(2q)}|a|^2e_{\I-z}h_{\I-z}^{1/(2q)})$.

We pause to justify this fact. Since this
justification parallels a similar justification in the proof of Proposition
\ref{cw}, we will be a little more terse here. As was the case in the proof of
Proposition \ref{cw}, $|a|^2$ is a uniform limit of a sequence of Riemann sums
of the form $\sum_{i=1}^n \lambda_i e_i$ with each $e_i$ a projection majorised
by $s_r(a)$. This uniform convergence in the compression $s_r(a)\M_1 s_r(a)$,
then ensures that the terms $h_1^{1/(2p)}(\sum_{i=1}^n \lambda_i
e_i)h^{1/(2p)}$ converge to $h^{1/(2p)}|a|^2h^{1/(2p)}$ in
$L^{2p}_{\varphi_1}(\M_1)$. Clearly $(e_z\sum_{i=1}^n \lambda_i e_i, e_{\I-z}\sum_{i=1}^n
\lambda_i e_i)$ will converge uniformly to $(|a|^2e_z, |a|^2e_{\I-z})$. Since this
convergence takes place in the compression $(s_r(a)e_z, s_r(a)e_{\I-z})(\M_1 \times
\M_1)(s_r(a)e_z, s_r(a)e_{\I-z})$, the terms
$$(h_z^{1/(2q)}e_z\sum_{i=1}^n \lambda_i e_ih_z^{1/(2q)},
h_{\I-z}^{1/(2q)}e_{\I-z}\sum_{i=1}^n \lambda_i e_ih_{\I-z}^{1/(2q)})$$will converge
to $(h_z^{1/(2q)}|a|^2e_zh_z^{1/(2q)}, h_{\I-z}^{1/(2q)}|a|^2e_{\I-z}h_{\I-z}^{1/(2q)})$

Thus by continuity there must exist a constant $M \geq 0$ so that
\begin{eqnarray*}
(\|[ah_z^{1/(2q)}]\|^{2q} + \|[ah_{(\I-z)}^{1/(2q)}]\|^{2q})^{1/q} &=&
\|(h_z^{1/(2q)}|a|^2h_z^{1/(2q)}, h_{(\I-z)}^{1/(2q)}|a|^2h_{(\I-z)}^{1/(2q)})\|\\
&\leq& M \|h_1^{1/(2p)}|a|^2h_1^{1/(2p)}\|\\
&=& M \|[ah_1^{1/(2p)}]\|^2.
\end{eqnarray*}

To conclude the proof we need only show that there exists some $K > 0$ with
$$\|[ah_J^{1/(2q)}]\| \leq K  (\|[ah_z^{1/(2q)}]\|^{2q} +
\|[ah_{(\I-z)}^{1/(2q)}]\|^{2q})^{1/(2q)}$$for all $a \in \M^{(0)}$, and apply
Lemma \ref{cwl2}. This fact is palpably clear if $q = \infty$, and hence we will
assume $1 \leq q < \infty$. For such a $q$ it is a simple matter to show that
$(r+s)^{1/q} \leq r^{1/q} + s^{1/q}$ for any $r, s \in \mathbb{R}_0^+$. Thus since
$h_z$ and $h_{\I-z}$ are commuting positive affiliated operators with $h_J = h_z + h_{\I-z}$,
it follows from the Borel functional calculus for such operators (see \cite[\S 5.6]{KR}),
that $h_J^{1/q} \leq h_z^{1/q} + h_{\I-z}^{1/q}$. Given any $a \in \M^{(0)}$, this in
turn has the effect of ensuring that $[h_J^{1/(2q)}a^*]^2 \leq [h_z^{1/(2q)}a^*]^2 +
[h_{\I-z}^{1/(2q)}a^*]^2$. Consequently $\|[h_J^{1/(2q)}a^*]^2\|_q \leq
\|[h_z^{1/(2q)}a^*]^2 + [h_{\I-z}^{1/(2q)}a^*]^2\|_q \leq \|[h_z^{1/(2q)}a^*]^2\|_q
+ \|[h_{\I-z}^{1/(2q)}a^*]^2\|_q$. Since $\|[h_J^{1/(2q)}a^*]^2\|_q =
\|[h_J^{1/(2q)}a^*]\|_{2q}^2 = \|[ah_J^{1/(2q)}]\|_{2q}^2$, and similarly
$\|[h_z^{1/(2q)}a^*]^2\|_q = \|[ah_z^{1/(2q)}]\|_{2q}^2$ and
$\|[h_{\I-z}^{1/(2q)}a^*]^2\|_q = \|[ah_{\I-z}^{1/(2q)}]\|_{2q}^2$, it is therefore
clear that
$$\|[ah_J^{1/(2q)}]\|_{2q} \leq (\|[ah_z^{1/(2q)}]\|^2 +
\|[ah_{(\I-z)}^{1/(2q)}]\|^{2})^{1/2}.$$

A simple application of H\"{o}lder's inequality then yields
\begin{eqnarray*} \|[ah_J^{1/(2q)}]\| &\leq&
(\|[ah_z^{1/(2q)}]\|^2 + \|[ah_{(\I-z)}^{1/(2q)}]\|^2)^{1/2}\\
&\leq& 2^{1/r}(\|[ah_z^{1/(2q)}]\|^{2q} + \|[ah_{(\I-z)}^{1/(2q)}]\|^{2q})^{1/(2q)}
\end{eqnarray*}
where $r \geq 1$ is chosen so that $\frac{1}{2} = \frac{1}{2q} + \frac{1}{r}$.

To see the final statement, observe that if $\mathcal{B} = J(\M_1)$, then $J$ induces
a Jordan $*$-isomorphism from $e\M_1$ onto $\mathcal{B}$. Thus in this case $e_z$
and $e_{\I-z}$ will indeed be disjoint and will respectively be mapped onto $z$ and
$J(\I)-z$ by $J$. Having centrally orthogonal supports, the weights $\varphi_z$ and
$\varphi_{\I-z}$ must commute.
\end{proof}

\end{document}